\documentclass[11pt]{amsart}
\usepackage{amsmath,amssymb,amsfonts,epsf,epsfig}
\usepackage[labelfont=normal]{subcaption}
\usepackage{enumerate,amsthm,color,mathdots}
\usepackage{euscript,graphicx}
\usepackage{url, float, hyperref, color} %, ulem}
\usepackage[usenames,dvipsnames,svgnames,table]{xcolor}
\usepackage{graphicx, mathtools}
\usepackage{tikz}

\newtheorem{theorem}{Theorem}[section]
\newtheorem{lemma}[theorem]{Lemma}
\newtheorem{proposition}[theorem]{Proposition}
\newtheorem{corollary}[theorem]{Corollary}

\newtheorem{definition}[theorem]{Definition}

\newtheorem{remark}[theorem]{Remark}

\renewcommand{\mod}{\hbox{{\rm mod}}\, }
\newcommand{\Aut}{{\rm Aut}}

\newcommand{\CC}{\mathbb{C}}
\newcommand{\SSS}{\mathbb{S}}
\newcommand{\gen}[1]{\langle #1 \rangle}

\newcommand{\GCD}{\mathrm{gcd}}
\newcommand{\beg}{\mathrm{beg}}
\newcommand{\inv}{\mathrm{inv}}

\newcommand{\Sym}{{\rm Sym}}

\newcommand{\C}{{\mathcal{C}}}

\newcommand{\ZZ}{\mathbb{Z}}

\newcommand{\A}{{\mathcal{A}}}
\newcommand{\B}{{\mathcal{B}}}

\newcommand{\D}{{\mathcal{D}}}
\newcommand{\Su}{{\mathcal{S}}}

\newcommand{\M}{{\mathcal{M}}}
\newcommand{\N}{{\mathcal{N}}}
\newcommand{\Pe}{\hbox{\rm P}}

\newcommand{\V}{{\rm V}}
\newcommand{\E}{\hbox{\rm E}}

\title[Cornerations in maps]{Transitive Cornerations in maps}

\author{Micael Toledo}
\address{Micael Toledo, Universit\'e Libre de Bruxelles, D\'{e}partement de Math\'{e}matique, C.P.216 - Alg\`ebre et Combinatoire, Boulevard du Triomphe, 1050 Brussels, Belgium}
\email{micaelalexitoledo@gmail.com}

\author{Alejandra Ramos-Rivera}
\address{Alejandra Ramos-Rivera, Institute of Mathematics, Physics and Mechanics, Jadranska 19, SI-1000 Ljubljana, Slovenia.}
\email{alejandra.rivera@upr.si}

\author{Primo\v z Poto\v cnik}
\address{Primo\v{z} Poto\v{c}nik, Faculty of Mathematics and Physics, University of Ljubljana, Jadranska 21, SI-1000 Ljubljana, Slovenia.\newline
\indent Also affiliated with: Institute of Mathematics, Physics and Mechanics, Jadranska 19, SI-1000 Ljubljana, Slovenia.
}
\email{primoz.potocnik@fmf.uni-lj.si}

\author{Stephen E. Wilson}
\address{Stephen E. Wilson, Northern Arizona University, Flagstaff, Arizona, USA}
\email{stephen.wilson@nau.edu}
\thanks{The research was partially supported by the Slovenian Research Agency, programme no. P1-0294 and research projects J1-1691 and J1-4351. The first author is thankful for financial support from the F\'ed\'eration Wallonie-Bruxelles -- Actions de Recherche Concert\'ees (ARC Advanced grant).}

\begin{document}

\maketitle

\begin{abstract}
A corner in a map is an edge-vertex-edge triple consisting of two distinct edges incident to the same vertex. A corneration is a set of corners that covers every arc of the map exactly once. Cornerations in a dart-transitive map generalize the notion of a cycle structure \cite{CS} in a symmetric graph. In this paper, we study the cornerations (and associated structures) that are preserved by a vertex-transitive group of automorphisms of the map.
\end{abstract}

%%%
\section{Introduction}
\label{sec:intro}

\color{Black}
Decompositions of complex structures into smaller and simpler substructures is a classical theme in mathematics. One instance
of this phenomenon in graph theory is the problem of decomposing graphs into subgraphs of specific type and in particular,
decompositions of graphs into cycles, where a number of difficult and deep results were proved and conjectured (consider, for
example the famous Oberwolfach problem \cite{Ober}).

Decomposition problems receive a somewhat different flavour when symmetry conditions are imposed.
In graph theory, a considerable attention received the study of decomposition (in fact, factorisations) of
graphs that are preserved under some vertex-transitive groups of automorphisms (see, for example, \cite{Hom1,Hom2,Muz}).

In \cite{CS}, a notion of {\em arc-transitive cycle decompositions} of graphs was introduced. In short, a cycle decomposition of a graph $\Gamma$
is a partition of its edge-set $\E(\Gamma)$ such that each part of the partition induces a connected 2-valent subgraph of $\Gamma$ (in short, a {\em cycle} of $\Gamma$).
Such a decomposition is {\em arc-transitive} provided that it is preserved by some group of automorphisms of $\Gamma$ acting transitively
on the {\em arc-set} (or the {\em dart-set}, as we call it in this paper), where an {\em arc} (or {\em dart}) is an unordered pair consisting of an edge and one of its end-vertices. (We may also think of an arc as a directed edge.)
For reasons of brevity we shall call an arc-transitive cycle decomposition a {\em cycle structure}.
%Cycle structures proved to be a very important tool in understanding {\em cubic vertex-transitive graphs} and {\em semisymmetric tetravalent graphs}.

The investigation of cycle structures began in \cite{CS}, where they were introduced as a means to understand 
an important class of tetravalent arc-transitive graphs, namely, those where the stabiliser of a vertex in the automorphism groups
acts on the neighbourhood as the dihedral group $D_4$.

 The results proved there were used in constructions of {\em semisymmetric tetravalent graphs} \cite{LRComb}  and played an instrumental role in the study of cubic vertex-transitive graphs \cite{CubicCensus}; the relationship with the cubic vertex-transitive graphs was described  in \cite{CubicCensus} via a {\em split graph} construction, which takes
 cycle structures in tetravalent graph as an input and returns a cubic vertex-transitive graph.

In the present paper, we broaden the scope in several directions.
While we keep the focus on finite graphs (and indeed tacitly assume the finiteness condition through the paper), 
we allow the graphs to have multiple edges (but no loops or semiedges). Further, we allow graphs to have higher valence. Next, we allow the parts of a cycle decomposition to consist of {\em circuits}, where, by a `circuit', we mean the (undirected and unrooted) trail of a closed walk in a graph that passes through each of its edges at most once (but perhaps several times through the same vertex). 

In fact, we shift our focus from `circuits' and `circuit decompositions' to `corners' and `cornerations'.
 Here, by a {\em corner} we mean an incident edge-vertex-edge  triple,
  and by a {\em corneration}, a set of corners that covers every incident vertex-edge pair of the graph exactly once.

Surprisingly, this simple change in perspective  simplifies many aspects of the investigation and constructions.  
  Formal definitions will appear in Section~\ref{sec:corn}.  The connections 
between these cornerations and circuit decomposition will be explained in 
Section~\ref{sec:cornerations} and explicitly codified in Theorem \ref{pro:C(L)}.
While corners and cornerations can be defined and studied in any graph, they 
become of particular interest when considered in the setting of (skeletons of) maps.  
In particular, a corner of a map divides the circular order of faces meeting that vertex 
into two set of consecutive parts.  If the valence is $q$ and the two parts have sizes 
$j$ and $q-j$, with $j\leq q-j$, we call it a $j$-corner, and then draw conclusions 
about the corner or  corneration to which it belongs relative to this $j$.

%%*** Primoz thinks that we should say something about special corners and cornerations, such as wedges and j-cornerations, here ***
 
 Throughout the paper, we will be interested in questions of symmetry.  If $L$ is a 
corneration in a map $\M$, let $\Aut(L)$ be the group of all symmetries of $\M$ 
which preserve $L$. If $\Aut(L)$ is transitive on $L$, we will call $L$ {\em transitive}. 

We begin to understand transitive cornerations by first examining their {\em local} 
behaviour, i.e. how corners around a vertex are arranged (the {\em local corneration}) and how the stabilizer of 
that vertex acts on them (the {\em local action}).  Our main results are Theorem \ref{th:locgps}, which says that a transitive corneration has one of three local actions, isomorphic to a $D_\frac{q}{2}, D_\frac{q}{4}$, or $C_\frac{q}{2}$, and Theorem \ref{pr:loccor}, which says that the local corneration is equivalent to one of two standard ones, determined by the parity of $j$.

In Sections~\ref{sec:wedge} and \ref{sec:STG}, we shift our focus to $1$-cornerations (or {\em wedge-cornerations} as we also call them here), as in some sense every $j$-corneration can be viewed
as a $1$-corneration in a related $j$-hole map (see Section~\ref{sec:Hj} for the definition of the $j$-hole operator and Remark~\ref{rem:jvswedge} for further discussion).
The main results about $1$-cornerations are Theorem \ref{th:faces}, which classifies the four ways that elements of a 1-corneration can appear in faces of the map, and Theorem \ref{the:class}, which uses the concept of symmetry type graphs to classify all transitive 1-cornerations of maps into 12 types.
 
 Section~\ref{sec:SymCor} is devoted to the study of {\em symmetric cornerations} of maps, that is, cornerations whose symmetry groups act
 transitively on the darts of the map. These cornerations are of particular interest as in the simple maps (that is maps whose skeleta are simple graph) the circuit decompositions that the cornerations induce are in fact arc-transitive cycle decompositions, as introduced in \cite{CS}.
In Theorem~\ref{the:symodd}, we prove that
a map $\M$ admits a symmetric $j$-corneration if and only if $j$ is odd and either $\M$ or the Petrie dual $\Pe(\M)$ 
admit a group of symmetries that acts half-reflexibly on $\M$ (see Section~\ref{sec:MapOp} for the definition of the Petrie dual and Section~\ref{sec:symflag} for the definition
of half-reflexible actions).

Finally, in Section~\ref{sec:SG}, we address one of our main motivations for this investigation, namely new constructions of vertex-transitive graphs of low valence.
We do this by introducing several versions of the {\em split graph constructions}, which generalise the split graphs as defined in \cite{CubicCensus}.
These constructions, in conjunction with existing censi of highly symmetrical maps, will be used in the forthcoming work on enumerating and cataloguing 
vertex-transitive graphs.

\section{Basic Definitions}\label{Basic}

\subsection{Conventions}

We will use exponential notation for permutations: if $\sigma$ is an element of a permutation group $G$ acting on a set $X$ and $x\in X$, we write $x^\sigma$ for the image of $x$ under $\sigma$, and $x^G$ for the orbit of $x$ under $G$.

For a positive integer $n$ let $\ZZ_n = \{0,1, \ldots, n-1\}$ be the ring of integers modulo $n$ with the convention that an arbitrary integer $a$ can be seen as an element of 
$\ZZ_n$ by identifying it with the remainder $a ~ \mod n$.

In many situations it will be convenient to adopt the notation $|a|_n$, read, ``the order of $a$ mod $n$'', to mean the additive order of the element $a$ in $\ZZ_n$.  It is equal
to $\frac{n}{\gcd(n,a)}$.

\subsection{Graphs}
In this paper, we will use the word {\em graph} to indicate what is often called a {\em multigraph}.  That is, we will allow graphs in this paper to have multiple (parallel) edges but no loops or semi-edges. Combinatorially, a graph will be
a triple $\Gamma:=(V,E,\partial)$ with $V$ and $E$ being disjoint sets of {\em vertices} and {\em edges}, respectively, and
$\partial\colon E\to {V\choose 2}$ being a function assigning a pair of distinct {\em end-vertices} to each edge. If $\Gamma$ is a graph $(V,E,\partial)$, then we let $\V(\Gamma)=V$ and $\E(\Gamma)=E$.
Distinct edges that have the same $\partial$-images are said to be parallel and graphs without pairs of  parallel edges are said to be {\em simple}.  \color{Black}   We define the {\em valence} of  a vertex $v$ of $\Gamma$ to be the number of edges $e$ such that $v\in \partial(e)$.  \color{Black}

It will sometimes be useful to consider a 
{\em topological form} $\Gamma^T$ of  $\Gamma$, which is obtained by taking for each edge of $\Gamma$, a space homeomorphic to a closed line segment, then labelling
the endpoints of each segment with the end-vertices of the corresponding edge, and finally identifying all endpoints having the same label.

A permutation $\alpha$ of $V\cup E$ preserving each of the sets $V$ and $E$ and satisfying $\partial(e^\alpha) =\partial(e)^\alpha$
for every $e\in E$ is called an {\em automorphism}  or {\em symmetry} of $\Gamma$. The set of all automorphism of $\Gamma$ forms
a subgroup of the symmetric group $\Sym(V\cup E)$, called the {\em automorphism group of $\Gamma$} and denoted $\Aut(\Gamma)$. If a subgroup $G$ of $\Aut(\Gamma)
$ acts transitively on $V$ or on $E$, then $\Gamma$ is said to be $G$-vertex- or $G$-edge-transitive, respectively (with the prefix $G$ omitted when $G=\Aut(\Gamma)$). 
Note that in a simple graph $\Gamma$, the action of $\Aut(\Gamma)$ on $\V(\Gamma)$ is
faithful and thus an automorphism group can be viewed as a permutation group on $\V(\Gamma)$.

A {\em dart} of a graph $\Gamma$ is a pair $\{v,e\}$ with $e\in \E(\Gamma)$ and $v\in \partial(e)$. If a subgroup
$G$ of $\Aut(\Gamma)$ in its induced action on the set $\D(\Gamma)$ of all darts of $\Gamma$ acts transitively, 
then $\Gamma$ is said to be {\em $G$-dart-transitive} (or simply {\em dart-transitive} if we want to avoid the reference to $G$).

A {\em walk} of length $k$ in a graph $\Gamma$ is a sequence $[v_0, e_0, v_1, e_1, v_2,\dots, e_{k-2}, v_{k-1}, e_{k-1}, v_k]$
with $v_i \in \V(\Gamma)$ and $e_i\in \E(\Gamma)$ such that $\partial(e_i) = \{v_i,v_{i+1}\}$ for all $i \in \{0, \ldots, k-1\}$.
If $v_k = v_0$ then we say that the above walk is {\em closed}. A closed walk of length $k$ is a  {\em circuit} provided that it has $k$ distinct edges and is a {\em cycle} provided that it has $k$ distinct vertices.

A {\em circuit decomposition} in a graph is a partition of its edges such that each element consists of the set of edges of a circuit.   A {\em circuit structure} is a circuit decomposition such that some dart-transitive group preserves the decomposition.  Similar definitions hold for {\em cycle decomposition} and {\em cycle structure}

\color{Black}

%%%
\subsection{Maps on surfaces}
\label{subsec:maps}

A {\em surface} $\Su$ is a connected compact 2-dimensional manifold (without boundary).  A \textit{map}  $\M$ is an embedding of $\Gamma^T$  for some graph  
$\Gamma$, called the {\em skeleton} of $\M$, into a surface $\Su$ so that each connected component of $\Su\backslash\Gamma^T$ is, 
topologically, a disk; these components are then called the  \textit{faces} of the map.  Henceforward, we will not distinguish between the graph and its topological version.
 If the  surface $\Su$ is orientable, we call the map \textit{orientable}.

The notions of vertices, valences, edges, darts, walks, cycles and circuits of a map are inherited from its skeleton.
\color{Black}     If $\Gamma$ is bipartite, we say that $\M$ is vertex-bipartite.
  If there is a partition of the faces of map  into two classes or 'colours' such that each face of a given colour border only those of the other, we say that the map is {\em face-bipartite}.  \color{Black}
Notice that the compactness of $\Su$ implies that the map has a finite number of vertices, faces and edges.

\subsection{Symmetries and operators on maps}
\label{sec:MapOp}
In discussing any phenomenon concerning maps, it will be useful to have 'operators', i.e., functions which relate one map to another, such as the classic {\em dual} operator.  Another is the Petrie operator: 
 \color{black}

%{\bf Petrie:} 

	A {\em Petrie path} in a map is a circuit in $\M$ such that for any three consecutive edges $e_1, e_2, e_3$, the edges $e_1$ and $ e_2$ lie in a face on one side of $e_2$, while $e_2$ and $e_3$ lie in the face on the other side of $e_2$.  Removing the interior of each face, and spanning each Petrie path with a new membrane, gives a new map, $\Pe(\M)$.  The faces of $\Pe(\M)$ are the Petrie paths of $\M$ and vice versa.

	 While the surface underlying $\Pe(\M)$ is usually not the same as the one underlying $\M$ itself,  the skeletons are the same.    If $\M$ is orientable then its Petrie paths must have even length.
   Of course, 
  $\Pe(\Pe(\M)) \cong \M$.

 \color{Black}
 We will see a third operator, called a $j$-hole operator and denoted by $H_j$, in the next section.

  \color{black}

%\subsection{Symmetry of maps}

  A \textit{symmetry} of a  map $\M$ is a symmetry of its skeleton induced by a homeomorphism of the underlying surface.   The group of all symmetries of $\M$ is denoted $\Aut(\M)$. 
If $\Aut(\M)$ acts transitively on the vertices of $\M$, we say that $\M$ is a vertex-transitive map.  
Further symmetry types of maps and groups acting on maps (such as reflexible, face reflexible and half reflexible) will be
  discussed in Section~\ref{sec:STG}. For general background on symmetry and operators of maps, we refer the reader to \cite{Wop}.

\subsection{Flags}
\label{sec:flags}
We  introduce the {\em barycentric subdivision} of a map: choose a point in the interior of each face to be its {\em center} and a point in the relative interior of each edge to be its {\em midpoint}.   Draw dotted lines to connect each face-center with each incidence of the surrounding vertices and edge-midpoints.  The original edges and these dotted lines divide the surface into triangles called {\em flags}.  Let $\Omega$ be the set of all flags in the map.  \color{Black}  

\begin{figure}[h]
\begin{center}
\includegraphics[width=0.50\textwidth]{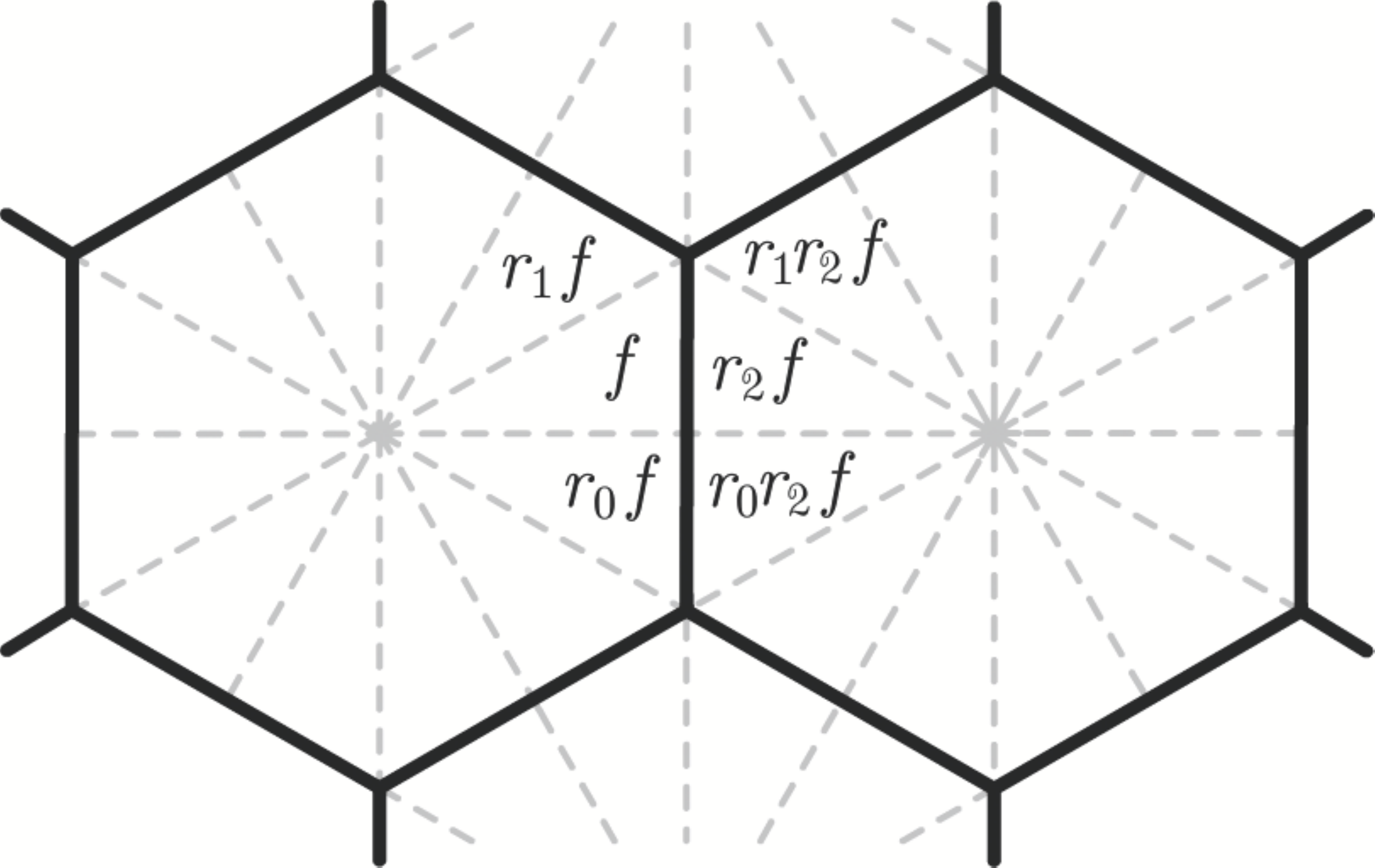}
\end{center}
\caption{Flags in a trivalent map with hexagonal faces.}
\label{fig:flaags}
\end{figure}

Note that a flag $f$ shares a side with three distinct flags (see Figure~\ref{fig:flaags}).
 One of them shares a face-center-to-edge-midpoint side with $f$; we will call this one
$r_0f$. One of the three shares a face-centre-to-vertex side with $f$; we will call this one $r_1f$. The third shares a vertex-to-edge-mid-point side with $f$; this one is denoted $r_2f$. In general, we shall call the flags $f$ and $r_if$ {\em $i$-adjacent} (note that this is well defined as $r_i(r_if) = f$). Moreover, note that $r_0$, $r_1$ and $r_2$ can be viewed as involutory permutations on $\Omega$ and that $r_0r_2 = r_2r_0$. Finally, connectedness of the map implies that the group $\langle r_0, r_1,r_2\rangle$ is transitive on $\Omega$.

\color{Black}

\subsection{Face- and half-reflexible maps}
\label{sec:symflag}
%\color{Black}Let $f$ be a flag in the map $\M$.  The picture in Figure \ref{fig:axes} shows flags near $f$, with certain axes in dashed lines, and numbered 0-0, 1-1, etc.\begin{figure}[h]\begin{center}\includegraphics[height=30mm]{FlagSym.pdf}\end{center}\caption{Axes in  a map}\label{fig:axes}\end{figure}

\color{Black}
Clearly, every symmetry of a map $\M$ induces a permutation on the set of flags $\Omega$ of $\M$. 
In fact, it can be seen that the induced action of $\Aut(\M)$ is faithful and semiregular (that is, a symmetry fixing one flag fixes every flag of the map).
In particular, a symmetry is determined by its action on a single flag. 

While there are many kinds of symmetry that a map may have, we wish to focus on just two.   First, we call $\M$ {\em reflexible} provided that $\Aut(\M)$ acts transitively on flags.   Second, we will call a subgroup $G$ of $\Aut(\M)$ {\em half-reflexible} provided that it is not transitive on flags, but for each face $F$ of $\M$, the stabilizer $G_F$ is transitive on the flags of $F$.   If $\Aut(\M)$ is half-reflexible we call $\M$ {\em half-reflexible}.  We will say a map $\M$ is {\em face-reflexible} if $\Aut(\M)$ has a half-reflexible subgroup.  From these definitions, the following facts   should be clear:

\begin{lemma}\label{lem:HFR}
For a face-reflexible 
map $\M$ and the half-reflexible subgroup $G$ of $\Aut(\M)$, each of the following holds: 

\begin{itemize}
%\item For any face $F$, $G_F$ is transitive on the edges of $F$.
\item   For any face $F$, all of its neighbors are in the same 
orbit under $G$ and that orbit does not contain $F$.
\item $\M$ is face-bipartite, and each partite set is an orbit under $G$. 
\item The index of $G$ in $\Aut(\M)$ is at most $2$.
\item If $G \neq \Aut(\M)$, then $\M$ is reflexible.
\end{itemize}
\end{lemma}

\begin{lemma}
\label{lem:HR}
A reflexible map $\M$ is face-bipartite if and only if it is face-reflexible.
\end{lemma}

\begin{proof}
From Lemma \ref{lem:HFR} we have that if $\M$ is face-reflexible, it is face-bipartite.
To complete the proof of this lemma, let us assume that $\M$ is reflexible and face bipartite.
Consider the faces in one class to be green, those of the other, red.  
Then every symmetry either sends red faces to red faces and greens to greens or it send each red face to a green  and every green to a red.   Let $G$ be the subgroup of 
$\Aut(\M)$ consisting of all symmetries of the first kind.   Then because each face is monochromatic, for each face $F$, $G_F$ acts transitively on the flags in that face.  Thus $G$ is half-reflexible and so $\M$ is face-reflexible..
\end{proof}

\color{Black}

%%%
\section{Corners in graphs and maps}
\label{sec:corn}
A set $\{e_1, v, e_2\}$ where $e_1$ and $e_2$ are distinct edges of a graph $\Gamma$ such that $v\in \partial_\Gamma(e_1) \cap \partial_\Gamma(e_2)$ is called a {\em corner} of $\Gamma$. 
If $g$ is a symmetry of a graph $\Gamma$ and $c=\{e_1,v,e_2\}$ is a corner, then
clearly $c^g:=\{e_1^g, v^g,e_2^g\}$ is also a corner of $\Gamma$. In this sense, $\Aut(\Gamma)$
acts upon the set of all corners of $\Gamma$.

\begin{figure}[hhh]
\begin{center}
\includegraphics[width=0.70\textwidth]{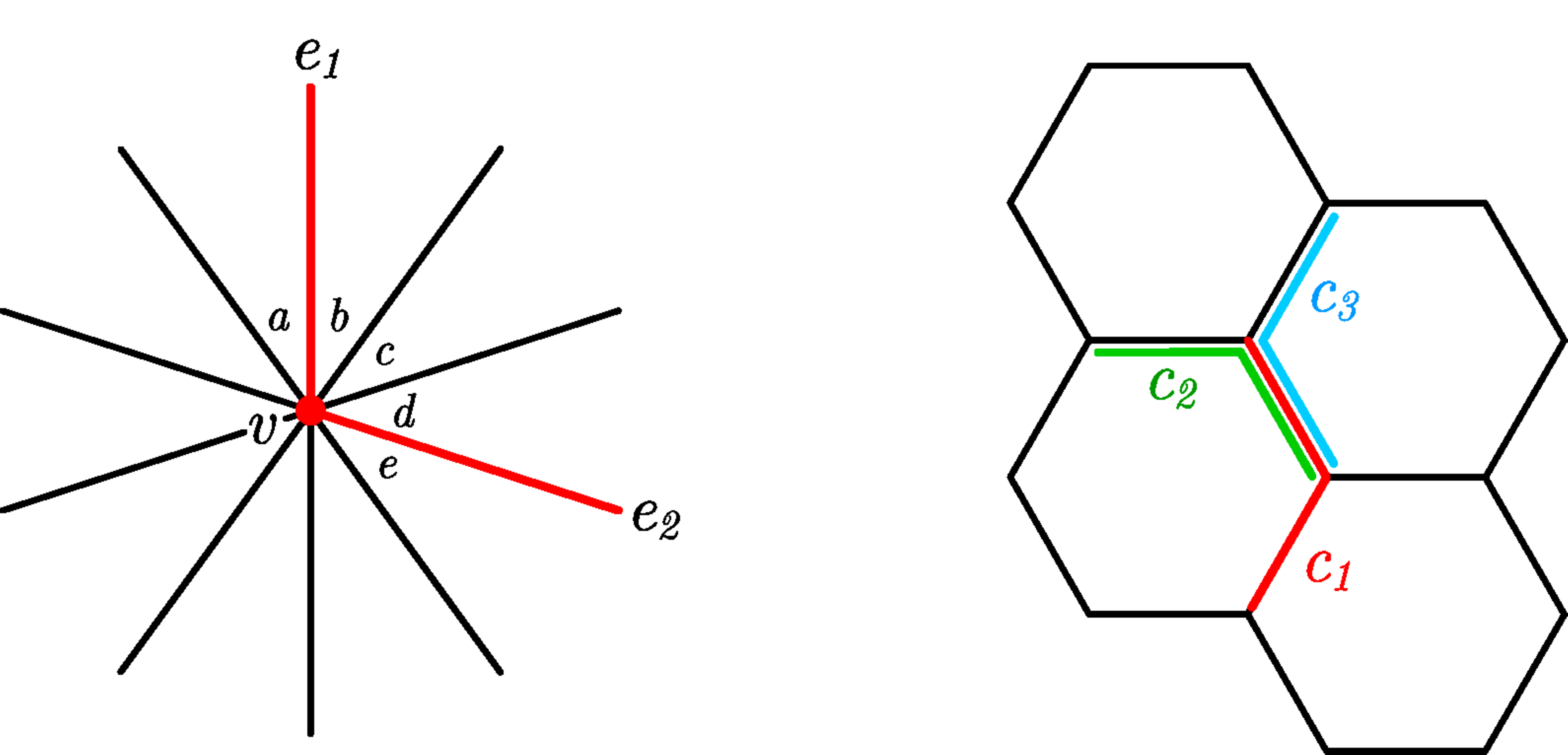}
\caption{On the left, a 3-corner $\{e_1,v,e_2\}$ in a map of valence 10. The wedges $b$, $c$ and $d$ are interior to $\{e_1,v,e_2\}$. On the right, three corners in a trivalent map. The corners $c_1$ and $c_2$ are convexly aligned; $c_1$ and $c_3$ are aligned with inflexion.}
\label{fig:CorPt}
\end{center}
\end{figure}

A {\em corner of a map} is defined to be a corner of its skeleton and a corner of a map $\M$ whose edges are consecutive in some face of $\M$  is called a {\em wedge} of $\M$.
If $v$ is a vertex of a map $\M$ of  valence $q$, then a corner $\{e_1, v, e_2\}$ separates some $j$, $1\le j \le \frac{q}{2}$, consecutive wedges at $v$ from the remaining $q-j$ wedges, which are also consecutive.
Such a corner is then called a $j$-corner. If $q$ is even, then a $\frac{q}{2}$-corner is said to be {\em straight}.
If $c$ is a $j$-corner which is not straight, then the set of 
$j$ consecutive wedges separated by $c$ from the rest of the wedges at $v$
is called the {\em interior} of $c$. \color{Black}    Further, we say that the corner {\em subtends} $j$ wedges. \color{Black} Figure \ref{fig:CorPt} shows a 3-corner $\{e_1, v, e_2\}$ in a map of valence $10$.  The three wedges $ b, c, d$ are the interior of the corner.
For a corner  $\{e_1, v, e_2\}$, the two wedges containing  the dart $\{e_1, v\}$ are {\em boundary wedges} for the corner, as are the two containing  $\{v, e_2\}$; for example, in Figure \ref{fig:CorPt}, wedges $a$ and  $b$ are boundary wedges, as are $d$ and $e$.

If $c_1 = \{e_1,v,e_2\}$ and $c_2 = \{e_2,u,e_3\}$, $u\not = v$, are two $j$-corners sharing the edge $e_2$, we say that $c_1$ and $c_2$ are {\em convexly aligned} provided that the interior wedges of each containing $e_2$ are on the same side of $e_2$; they are {\em aligned with inflection} otherwise (see the right-hand side of Figure \ref{fig:CorPt}).

Given a finite sequence  $\{c_i = \{e_i, v_i, e_{i+1}\} | i \in \ZZ_r\}$ of $j$-corners, it might be that every pair of consecutive corners is convexly aligned, in which case, we call the sequence a {\em $j$-th order hole}.  If, instead, every consecutive pair is aligned with inflection, we say it is a {\em $j$-th order Petrie path}.  Then  $j$-th order holes generalize faces, in the sense that a face is a first-order hole; similarly, $j$-th order Petrie paths generalize Petrie paths.

\subsection{$j$-th hole operator}
\label{sec:Hj}
Let us conclude the section by defining the $j$-th hole operator $H_j$ on maps, first introduced in \cite{Wop}.  For any map $\M$, dissolve its faces, leaving only the skeleton of $\M$, and then span 
each circuit of the skeleton induced by a $j$-hole of $\M$ with a membrane.

If $d = \GCD(j,q)$ is equal to 1, then the result is a map on a surface, which we can call $H_j(\M)$.  If  $d>1$, then the resulting topological space is not a surface but rather surface pinched at every vertex of the skeleton. By resolving the singularities by separating each vertex into $d$ new vertices, each of valence $\frac{q}{d}$, one obtains a number of disjoint maps.  Assign  the name $H_j(\M)$ to the set of all the component maps.

We could make this definition more rigorous by defining maps combinatorially in terms of set of flags $\Omega$ and the involutions $r_0$, $r_1$ and $r_2$ as described in Section \ref{sec:flags}. Then $H_j$ is the set of maps corresponding to the orbits under the action of the group generated by $r'_0$, $r'_1$ and $r'_2$ where $r'_0 = r_0$, $r'_1 = r_1(r_2r_1)^{j-1}$ and $r'_2 = r_2$. 

For integers $j, k$, it should be clear that  $H_j(H_k(\M)) = H_{jk}(\M)$.  Thus if $j$ and $q$ are relatively prime and $\N = H_j(\M)$, then there is an integer $k$ such that $\M = H_k(\N)$.

\section{Cornerations}
\label{sec:cornerations}

Note that each corner $\{e_0,v_0,e_1\}$ contains two darts, namely, $\{e_0,v_0\}$ and $\{v_0,e_1\}$.
A {\em corneration} of a graph $\Gamma$
is a set $L$ of corners of $\Gamma$ such that every dart of $\Gamma$ is contained in exactly one corner of $L$; it follows that then every edge belongs to exactly two corners of $L$.
A {\em corneration of a map} is then defined as a corneration of its skeleton.

%% \primozcomment{What follows is essentially a claim that decomposition of edges in circuits is equivalent to cornerations}

If $\C$ is a circuit decomposition of a graph $\Gamma$, then a corner $\{e_0, v_0, e_1\}$ such that $e_0,v_0,e_1$
is a consecutive triple in some circuit of $\C$ is called a {\em $\C$-corner}. Clearly, \color{Black} the set $L(\C)$ of  \color{Black} all $\C$-corners is a corneration.

Suppose now that $L$ is a corneration of a graph $\Gamma$.
Given a dart $d_0 = \{e_0, v_0\}$, there is a unique edge $e_1$ such that $\{e_0, v_0, e_1\}\in L$.   Denote this  $\{e_0, v_0, e_1\}$ as $L(d_0)$.  Now,   $v_0$ is one endvertex of $e_1$; let $v_1$ be the other.  We define $d_1 = \{e_1, v_1\}$ to be the {\em predecessor} to $d_0$.  This predecessor to $d_0$ is unique and, similarly, each dart has a unique successor.  

This allows to define $d_{i+1}$ to be the predecessor of $d_i$ for every $i \geq 0$. Let $v_i$ and $e_i$ be such that $d_i = \{v_i,e_i\}$. Since the number of darts is finite there exists a smallest integer $k>0$ such that $d_k = d_0$. By uniqueness of successor and predecessor  we see that the sequence $(v_0,e_0,\ldots,e_{k-1},v_{k})$ is a circuit in $\M$. Call it $O(d_0)$. Note that $O(d_0) = O(d)$ of every dart $d$ contained in the circuit $O(d_0)$.
 %Beginning with $d_0$ and recursively defining $d_{i+1} $ to be the predecessor to $d_i$, we get a dart walk which must be closed by the finiteness of  $\M$. 
  %Uniqueness of successor and predecessor guarantees that the first repetition of an edge will happen when $d_k  = d_0$ for the first time.   Thus the edges in this walk must be a circuit in $\Gamma$, call it $O(d_0)$.  
  %Note then, that  if the endvertices of $e_0$ are $v_0$ and $x$, then $\{e_0, x\}$ is the successor to $\{e_1, v_0\}$.   Thus, the circuit  $O(d_0)$ is the same as $O(\{e_0, x\})$.    
  Let $\CC(L)$ be the set of all $O(d)$ such that $d$ is a dart of $\Gamma$.
We present, without proof, the following theorem,  which follows easily from the above definitions.

\begin{theorem}
\label{pro:C(L)}
If $L$ is a corneration in the graph $\Gamma$, then $\CC(L)$ is a circuit decomposition of $\Gamma$.  Conversely, if $\C$ is a circuit decomposition of $\Gamma$,
then  $L(\C)$ is a corneration $L$ in $\Gamma$ such that $\C = \CC(L)$. 
\end{theorem}

Let now $\M$ be a map of even valence $q$.
Then the set  $L$ of all straight corners of $\M$ is a corneration, and the circuits in $\CC(L)$ are called  {\em lines.}

If $L$ is any corneration of a map $\M$,
then we will call $L$ {\em uniform} provided that for some $j\leq \frac{q}{2}$, every element of $L$ is a $j$-corner. 
  We might also say that $L$ is $j${-\em uniform} or that $j$ is the {\em width} of $L$. For such $L$ with $j < \frac{q}{2}$, define its {\em $j$-complement} to be the set of all $j$-corners  {\em not} in $L$. Note that the $j$-complement of a $j$-corneration is then also a $j$-uniform corneration.  

\color{Black}
\begin{remark}\label{rem:comp}
If $L$ is $j$-uniform for some $j < \frac{q}{2}$, consider a corner $c_0 = \{e_0, v, e_1\}$ of $L$.  The other $j$-corner containing $\{v, e_1\}$ is some $c_1 = \{e_1, v, e_2\}$, which then must be in the $j$-complement of  $L$.  Proceeding inductively, we get a sequence of $k$ corners $c_i = \{e_i, v, e_{i+1}\}$ for $i =0,1,2,\dots,k-1$, where $c_{k-1} = \{e_{k-1}, v, e_0\}$.  In this sequence,  $c_i\in L$ holds exactly when $i$ is even.   It follows that $k$ is even. % As this is true for each $j$-corner at $v$, we see that $q$ must be even. 

 This sequence may or may not include all of the $j$-corners at $v$.  More precisely, the number $k$ is $|j|_q$.  Thus the sequence uses all $j$-corners at $v$ if and only if $j$ and $q$ are relatively prime.   In particular, if this happens, then $j$ is odd (as the valence $q$ is even in every graph admitting a corneration).
\end{remark}
\color{Black}
\begin{remark}
\label{rem:disconnect}
Let $\N = \Pe(\M)$.  If $c$ is a $j$-corner of $\M$ then it is also a $j$-corner of $\N$, and has the same interior wedges and boundary wedges in $\N$ as in $\M$.  If $L$ is a corneration of $\M$ it is also a corneration of $\N$.   Further, let $\N' = H_j(\M)$.  If $c$ is a $j$-corner of $\M$ then it is also a $1$-corner in $\N'$.  If $L$ is a  $j$-corneration of $\M$ it is also a 1-corneration of $\N'$. We note that if $j$ is not coprime to $q$, then $\N'$ might be a union of several maps. In this case the definition of a corneration can be extended in a obvious way to $\N'$.

\end{remark}

\section{Symmetry in Cornerations}\label{sec:cornsym}

If $L$ is a corneration of a graph $\Gamma$ (or a map $\M$), let $\Aut(L)$ be the subgroup of $\Aut(\Gamma)$  (or of $\Aut(\M)$)
 consisting of all symmetries of $\Gamma$ (or of $\M$, respectively) which preserve $L$. 
 
Let $G$ be a subgroup of $\Aut(L)$ which is transitive on corners of $L$.  
Then we shall say that $L$ is $G$-transitive.
It might happen that $G$ is transitive on the darts of $\Gamma$ (or on $\M$), in which case, we call $L$ 
{\em $G$-symmetric}. 
Since every $G$-transitive or $G$-symmetric corneration is also $\Aut(L)$-transitive or symmetric (respectively),
we may omit $G$ in the above definitions, and simply say that $L$ is transitive, or symmetric, respectively.
The first of the following two propositions is obvious while the second follows easily from Theorem~\ref{pro:C(L)}.

\begin{proposition}
\label{pr:uniform}
Every transitive corneration of a map is uniform.  
\end{proposition}

\begin{proposition}
If $L$ is a  corneration of a map $\M$, and  $\C = \CC(L)$, then $L$ is symmetric if and only if $\C$ is a circuit structure.  
If, in addition, $\M$ is simple then $\C$ is a cycle structure.
\end{proposition}

\begin{remark}
Following Remark \ref{rem:disconnect}, let $\N'=H_j(\M)$.   If $L$ is a transitive corneration of $\M$, then it is a transitive corneration of $\N'$.  However, if $\N'$ is disconnected, then $L$ might be transitive (or even symmetric) on $\N'$ and still not be transitive on $\M$.
\end{remark}

\begin{lemma}\label{lem:DTsym}
Let $\M$ be a 
%dart-transitive  
map of valence $q$. If there exists a symmetric $j$-corneration of $\M$ with $j < \frac{q}{2}$, then $j$ is odd.
%
%Let $\M$ be a 
%\color{Black} simple? \color{black}
 %map and let $L$ be a corneration of $\M$. If $L$ is symmetric and $j$-uniform for some  then $j$ is odd.
 \end{lemma}

\begin{proof}
Suppose that $L$ is a symmetric $j$-corneration of $\M$ and that $j < \frac{q}{2}$.
Let $c=\{e_1,v,e_2\}$ be a corner of $L$ and suppose, for a contradiction, that $j$ is even. Then there exists an edge $e_3$ incident to $v$ and interior to $c$ that is ``$\frac{j}{2}$ wedges apart'' from $e_1$ and from $e_2$. That is, the corner $\{e_1,v,e_3\}$ (which is not in $L$) subtends $\frac{j}{2}$ wedges, as does the corner $\{e_2,v,e_3\}$. Since $L$ is symmetric, there must be a $g \in \Aut(L)$ mapping the dart $\{e_1,v\}$ to $\{e_2,v\}$ and since it fixes $c$, it must map $\{e_2,v\}$ back to $\{e_1,v\}$. We see then that $g$ acts like a reflection swapping $e_1$ and $e_2$ while fixing $e_3$ and $v$. Let $e_4$ be the unique edge such that $c':=\{e_3,v,e_4\}$ is a corner of $L$ and observe that $g$ must fix $c'$ point-wise. It follows that $g$ swaps every wedge interior to $c'$ with an exterior wedge, since it swaps $e_1$ with $e_2$. This can only happen if the corner $c'$ subtends exactly $\frac{q}{2}$ wedges, contradicting the assumption that every element of $L$ has width $j< \frac{q}{2}$.
\end{proof}
\color{Black}

%\begin{remark}
%Observe that the converse of Lemma~\ref{lem:DTsym} does not hold in general. In particular, the $4$-valent toroidal map  $\{4,4\}_{2,1}$ (see \cite{CoxMos} for the meaning of this notation) with five vertices and five faces admits no symmetric $1$-cornerations.
%It would be interesting to investigate under which additional conditions on $\M$ and on $j$, the existence of a symmetric $j$-corneration can be guaranteed, and whether all such can be classified. While this exceeds our ambitions in this paper, we would like to propose it as a possible future avenue of research.
%\primozcomment{Change this remark. We will now state and prove partial converse of Lemma~\ref{lem:DTsym}}
%\end{remark}

The remainder of this paper is in two parts: first, we discuss the group of a transitive corneration and discuss how it acts on vertices and on faces of the map; second, we introduce the generalizations of the split graph and consider the connectedness of the constructed graphs.

\color{Black}
%%   This is not true   \begin{theorem}
%%Let $\M$ be a simple map and let $L$ be a corneration of $\M$. If $L$ is symmetric, then $L$ is $j$-uniform and one of the following occurs:
%\begin{enumerate}
%\item $j$ is odd, $j < \frac{q}{2}$ and $\M$ is face-bipartite
%\item $j = \frac{q}{2}$
%\end{enumerate}
%\end{theorem}

%\primozcomment{We need to prove that}

\section{Local Actions  and Local Cornerations}\label{sec:local}   \color{Black}

%\primozcomment{Reconsider the content of this section. There's probably some overlap with later section. Some parts might need to go to the split graph section}
\subsection{Local cornerations and local actions}

In classifying actions of groups on a map   one important step is to determine the ways in which the stabilizer of a vertex  acts on its neighbourhood. Loosely speaking, this is the {\em local action} of the group.  In this section, we work to determine the possible local actions of a  group which is transitive on a corneration.  \color{Black}

First, notation:   If $L$ is a corneration of a map $\M$ and $v$ is a vertex of $\M$, let $L_v$ be the set of corners in $L$ containing $v$

and call this set the {\em local corneration} at $v$.  Number the edges incident to $v$ with $0,1,2,\ldots, q-1$ consecutively around $v$, and consider those labels to be in $\ZZ_q$.

There are, of course, $q$ choices of where to start the numbering and two choices for which way to go around the circle.  The $2q$ distinct numberings  will be  considered to be {\em equivalent}.

Let $\rho$ be the permutation $i\rightarrow i+1$,  and for each $f\in \ZZ_q$, let $\mu_f$ be the permutation $i\rightarrow f-i$.    Then $\rho$ acts as a rotation about $v$ and each $\mu_f$ as a reflection.   Those permutations generate a group isomorphic to the dihedral group isomorphic to $D_q$, and they may or may not extend to symmetries of $\M$.

Now consider a transitive (and therefore uniform) $j$-corneration $L$ of $\M$ for some $j<\frac{q}{2}$.  Let $G=\Aut(L)$ and think of the vertex stabilizer $G_v$ as a permutation group acting on the labels $0,1,2,\ldots, q-1$ of edges incident to $v$. Then clearly $G_v \leq D_q$. It should now be clear that the following are true:

\begin{enumerate}
\item We can represent the corners at $v$  by pairs $c_i = (i, i+j)$ for $i \in \ZZ_q$.  There are $q$ such pairs and half of them are in $L_v$.
\item $G_v$ cannot contain $\rho$, because for any $c_i$ in $L_v$, $c_i^{\rho^j}$ is in $L_v$, and shares the edge $e_{i+j}$ with it.
\item If $f$ is even, say $f = 2f'$, then $G_v$ cannot contain $\mu_f$ because $\mu_f$ would fix the edge $f'$ but not the other edge in its corner.
\item Because $G_v$ is transitive on the elements of $L_v$, its order must be divisible by $\frac{q}{2}$ and $G_v$ has at most 2 orbits on the edges at $v$.
\item If $G_v$ has two orbits on edges, then each corner $c_i $ of $L_v$ must contain one edge from each orbit.
\end{enumerate}

There are only three subgroups of $D_q$ that satisfy these requirements:
\begin{enumerate}
\item   $D_\frac{q}{2}$, generated by $\mu_1$ and $\mu_3$, containing $\rho^2$; this subgroup is transitive. %on edges incident to $v$.
We say the local action is {\em half-dihedral}. (HD, for short)  
\item  $C_\frac{q}{2}$, generated by $\rho^2$; it has two orbits of edges, those that are even and those that are odd. 
 We say the local action is {\em half-cyclic} (HC, for short)  

\item  $D_\frac{q}{4}$,  generated by $\mu_1$ and $\mu_5$, containing $\rho^4$.  Here there are two orbits of edges; one consists of edges with numbers  which are 0 or 1 mod $4$, and one has those which are 2 or 3 mod $4$.
\color{Black}  We say the local action is {\em quarter-dihedral}.  (QD, for short) \color{Black}
\end{enumerate}

Since $\Aut(L)$ is vertex-transitive by definition, the type of local action is independent of the choice of $v$. Then following theorem summarizes the discussion above.

\begin{theorem}
\label{th:locgps}
If $L$ is a transitive corneration in a map $\M$ of valence $q$, and $G = \Aut(L)$, then its local action is of type {\rm HD},  {\rm HC} or  {\rm QD}.
As a consequence, the index of $G$ in $\Aut(\M)$ is $1$, $2$ or $4$. Further, if $L$ is symmetric, then the local action is of type {\rm HD}.
\end{theorem}

%\subsection{Local Cornerations}

It will be convenient to define two standard local cornerations:
\begin{definition}
Using, as before, $c_i$ to stand for the $j$-corner $(i, i+j)$  at vertex $v$ in map $\M$,  we define:
\begin{itemize}
\item For odd $j$ and even $q$, the {standard odd local corneration} is $\{c_{2i} \mid i\in \{0,\ldots,\frac{q}{2}-1\}\}$.
\item For $j\equiv 2$ (mod 4) and $q$ divisible by 4,  the {standard even local corneration} is $\{c_{i} \mid i  ~ { \rm mod ~4}  ~\in \{0,3\}  \}$.
\end{itemize}

\end{definition}

With these definitions established, the following theorem holds:
\begin{theorem}\label{pr:loccor}\color{Black} 
If $L$ is a transitive $j$-corneration of a map $\M$ and $v$ is a vertex of $\M$, then
\begin{itemize}
\item $j$ is odd and $L_v$ is equivalent to the standard odd local corneration,

or

\item $j\equiv 2$ (mod 4), $L_v$ is equivalent to the standard even local corneration and the local action of $\Aut(L)$ is of type $\textrm{QD}$.
\end{itemize}
\end{theorem}

\begin{proof}
Let $G = \Aut(L)$.  We consider local cornerations having each of the groups mentioned in Theorem \ref{th:locgps}.

\color{Black}

Suppose first that $G_v$ is isomorphic to $C_\frac{q}{2}$.  Then a corner in $L_v$ must contain one even and one odd edge; thus $j$ must be odd.  Assume with no loss of generality, that $(0, j)$ is in $L_v$.   Then $L_v$ must contain and consist of all corners of the form $(2i, 2i+j)$ for $i \in \ZZ_q$.   This is exactly the  standard odd local corneration.

Now suppose that $G$ is isomorphic to $D_\frac{q}{2}$.   Then $j$ must be odd, for otherwise $\rho^j$, which is in $\gen{\rho^2}$, would send each $(i,i+j)$ in $L_v$ to $ (i+j, i+2j)$ also in $L_v$, which would not satisfy the corneration definition.  Renumbering the edges if necessary, we can assume that $(0,j)\in L_v$, and then because $\rho^2\in G$, we have that $L_v$ is equivalent to the standard odd local corneration.

Finally, suppose that $G$ is isomorphic to $D_\frac{q}{4}$.  In particular, suppose that $G$ is generated by $\mu_1$ and $\mu_5$, and so contains $\rho^4$.  Then, as noted above, one orbit of edges is the set of those which are 0 or 1 (mod 4), and the other, those that are 2 or 3 (mod 4).  It follows that there are three orbits of wedges:
\begin{enumerate}
\item Those of the form $(4i, 4i+1)$
\item Those of the form $(4i+2, 4i+3)$
\item Those of form $(4i+1, 4i+2)$ or $(4i-1, 4i)$
\end{enumerate}

We can generalize this comment a little:  For $j\equiv 1$ (mod 4), there are three orbits of corners of width $j$:

\begin{enumerate}
\item Those of the form $c_{4i}$
\item Those of the form $c_{4i+2}$
\item Those of form $c_{4i+1}$ or $c_{4i-1}$
\end{enumerate}

Only the third of these covers all edges at $v$, and so only this can be a local corneration, but it is equivalent to the standard odd local corneration, with which we are familiar.

On the other odd hand, for $j\equiv 3$ (mod 4), there are three orbits of corners of width $j$:

\begin{enumerate}
\item Those of form  $c_{4i-1}$
\item Those of the form $c_{4i+1}$
\item Those of the form $c_{4i}$ or $c_{4i+2}$
\end{enumerate}
Only the third of these covers all edges at $v$, and so only this can be a local corneration.But it is identical to the standard odd local corneration, with which we are familiar.

Now, the width cannot be divisible by 4, because  edges which are equivalent mod 4 are in the same orbit.  Thus the only remaining possibility for $j$ is 2 (mod 4).   Suppose, without any loss of generality, that $(0,j)\in L_v$; then $L_v$ also contains $(0,j)^{\mu_{j-1}} = c_{-1}$ and thus $L_v$ contains all pairs $(i, i+j)$ such that $i\equiv 0,3$ (mod 4).  
Thus, corners in $L_v$ come in consecutive pairs: $c_{4i-1} = (4i-1, 4i-1+j)$ and $c_{4i} = (4i, 4i+j)$, while the corners  $(4i+1, 4i+1+j)$ and $(4i+2, 4i+2+j)$ are not in $L_v$.  Thus $L_v$ is equivalent to the  standard even local corneration.

\end {proof}
\color{Black}

\section{Wedge Cornerations}
\label{sec:wedge}

In this section, we consider cornerations in which each corner is a 1-corner, also called a  wedge. This restriction is not as severe as it might appear, as Remark \ref{rem:comp} allows us to analyse a $j$-corneration of $\M$ by examining wedge cornerations of $H_j(\M)$ where $j$ is relatively prime to $q$. Note that each wedge {\em belongs} to a unique face in the sense that the two edges of the wedge appear as two consecutive edges on the boundary of a unique face.

 Suppose that $L$ is a  corneration of  a map $\M$ whose elements are wedges, and let $G=\Aut(L)$.    For brevity, we refer to wedges in $L$ as {\em in-wedges} and those not in $L$ as {\em out-wedges}.

\color{Black}
  Given a face $f$ of length $p$, \color{Black} here are some ways that in-wedges might appear around $f$:
\begin{itemize}
\item{(A)} Every wedge of $f$ is an in-wedge.
\item{(B)} $p$ is even and every other wedge of $f$ is an in-wedge.
\item{(C)} $p$ is divisible by 3 and every third wedge of $f$ is an out-wedge.
\item{(D)} $p$ is divisible by 4 and pairs of in- and out-wedges alternate.
\item{(E)} Every wedge of $f$ is an out-wedge.
\end{itemize}
We say that $f$ has {\em pattern} A, or B, etc..  These possibilities are illustrated in Figure \ref{fig:XABCD}, where the in-wedges are indicated by dots.

\begin{figure}[hhh]
\begin{subfigure}{2in}
\begin{center}
\includegraphics[height=19mm]{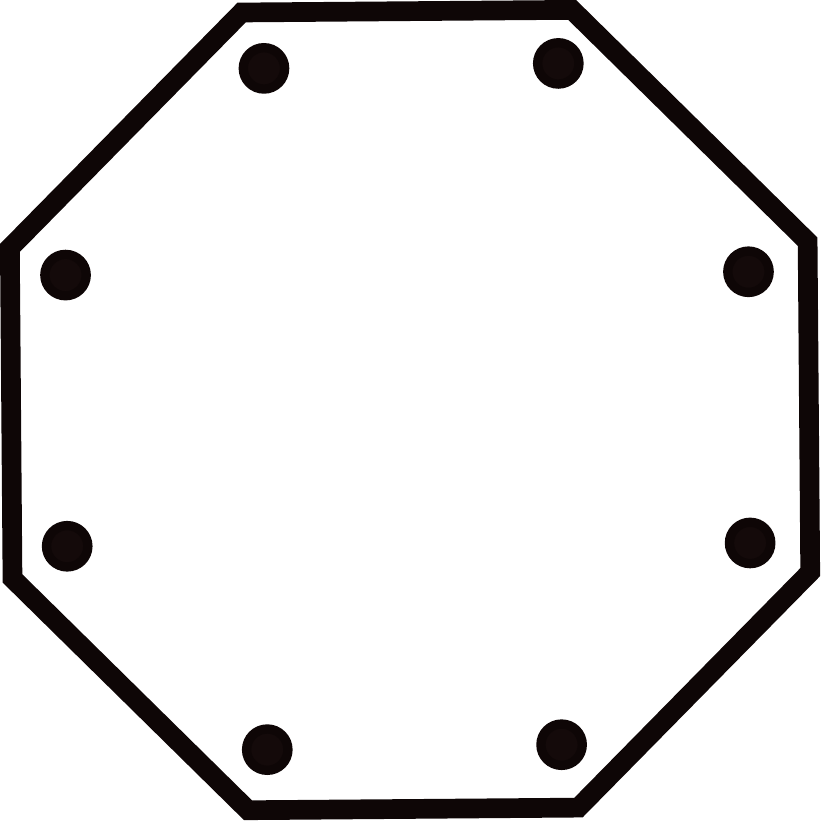}
\caption{Pattern A}\label{fig:XA}
\end{center}
\end{subfigure}
~
\begin{subfigure}{2in}
\begin{center}
\includegraphics[height=19mm]{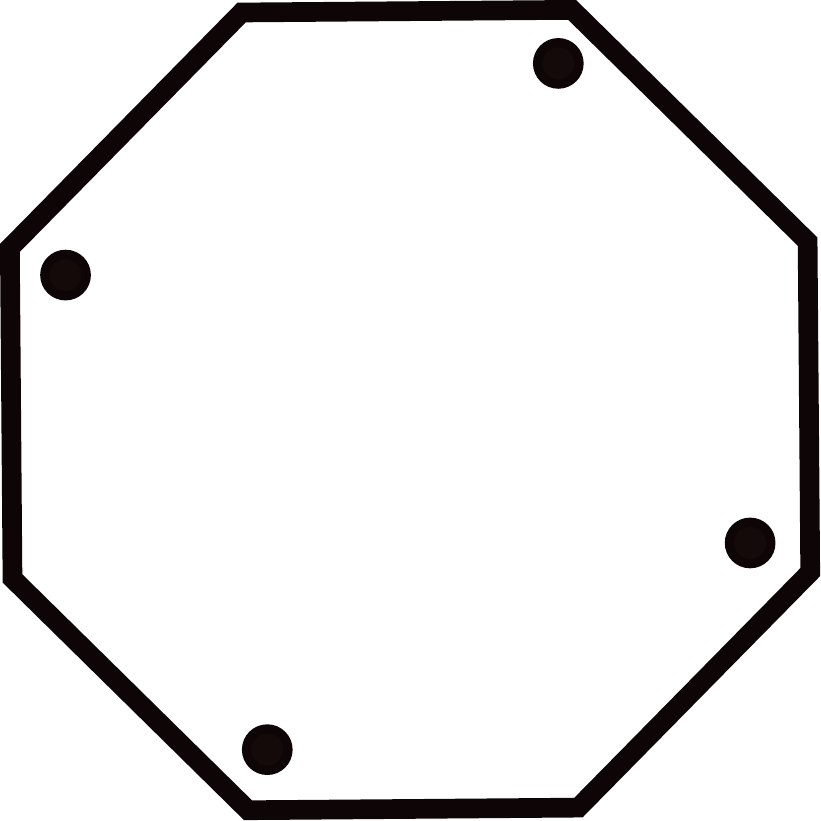}
\caption{Pattern B}\label{fig:XB}
\end{center}
\end{subfigure}
 ~
\begin{subfigure}{2in}
\begin{center}
\includegraphics[height=19mm]{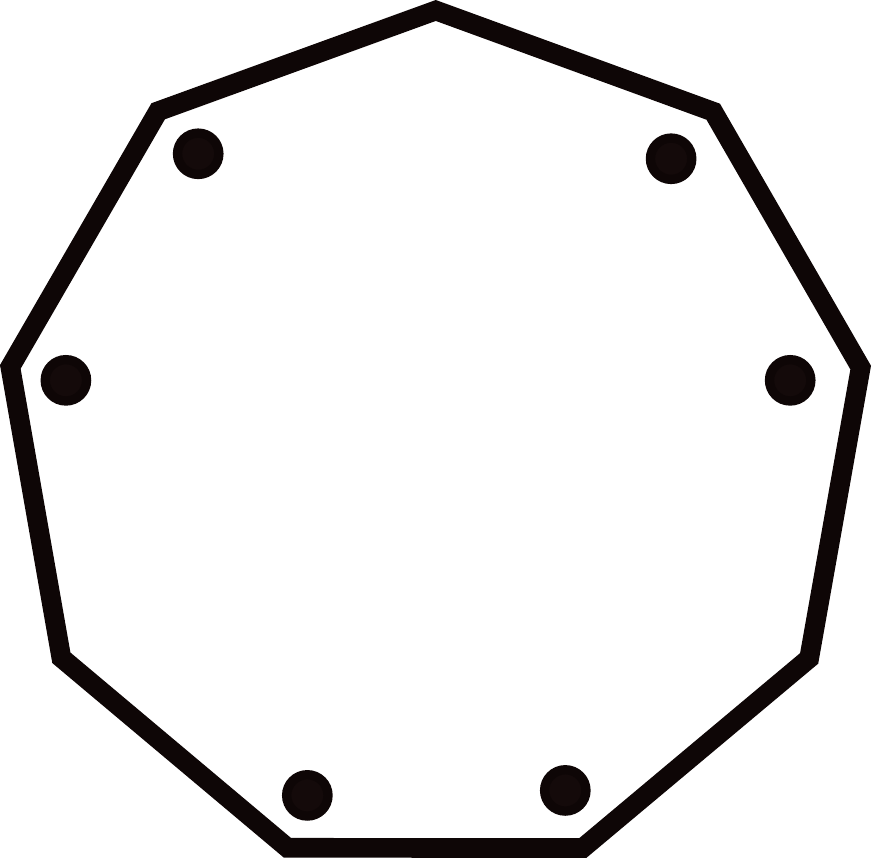}
\caption{Pattern C}\label{fig:XC}
\end{center}
\end{subfigure}

\begin{subfigure}{2in}
\begin{center}
\includegraphics[height=19mm]{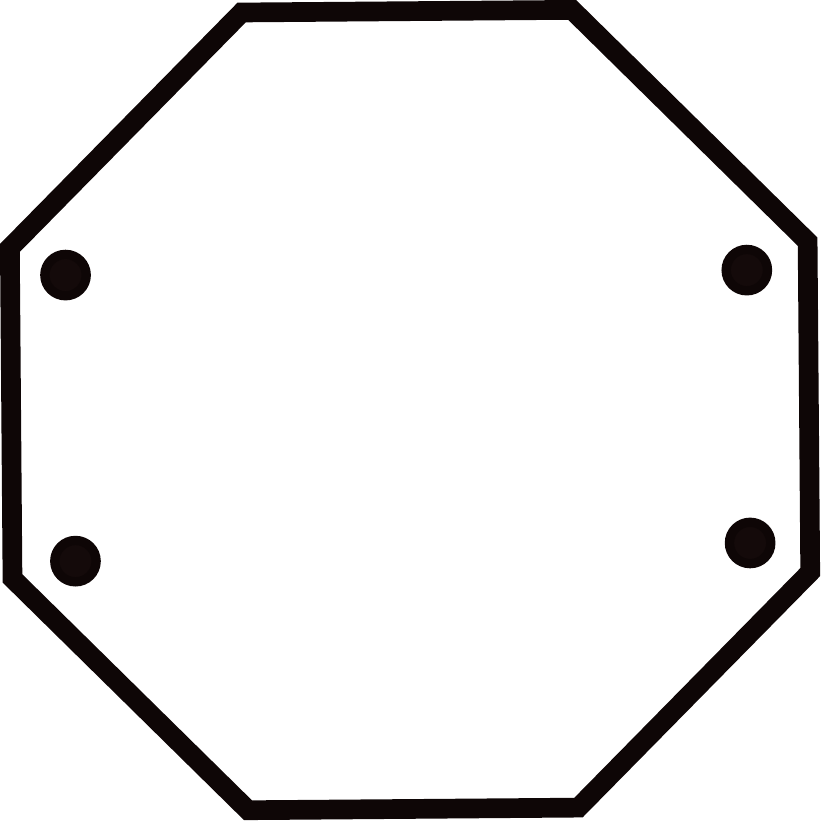}
\caption{Pattern D}\label{fig:XD}
\end{center}
\end{subfigure}
~
\begin{subfigure}{2in}
\begin{center}
\includegraphics[height=19mm]{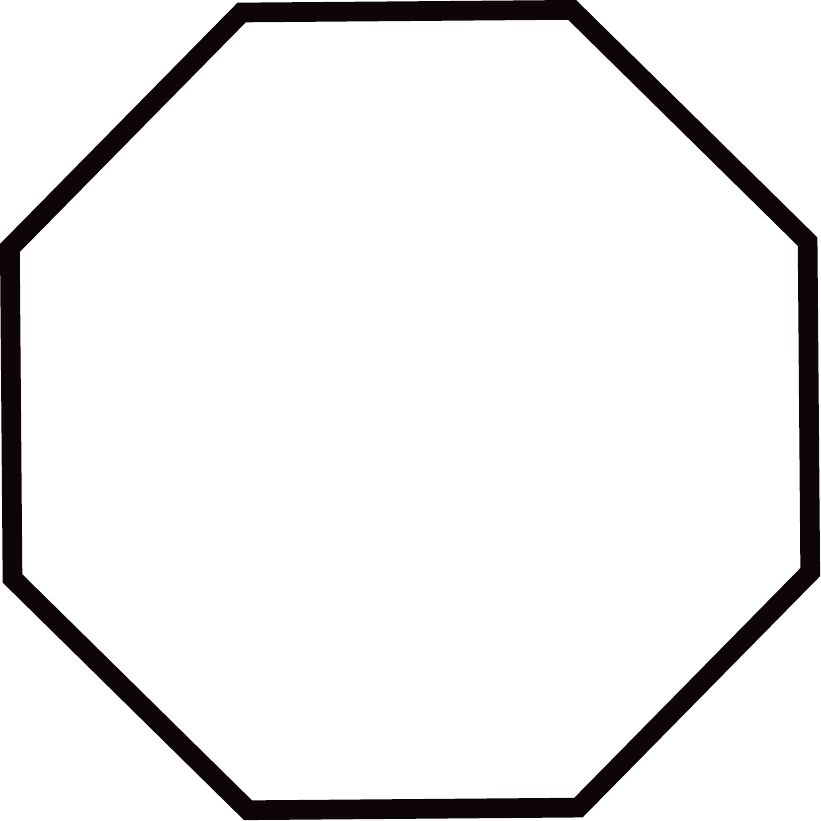}
\caption{Pattern E}\label{fig:XE}
\end{center}
\end{subfigure}

\label{fig:XABCD}
\caption{Patterns of corners in a face}\label{fig:XABCD}
\end{figure}

\begin{theorem}\label{th:faces}
Given that $L$ is a transitive corneration of  a map $\M$  whose elements are wedges there are four possibilities for the way that patterns of in-wedges appear in faces of $\M$:
\begin{enumerate}
\item $\M$ is face-bipartite and all faces of one color have pattern A, while the rest have pattern E.
\item Every face has pattern B.
\item Faces are of two patterns, C and E. Faces with patern E border only those with pattern C, while a face with pattern C meets faces with patterns C and E.
\item Every face has pattern D.
\end{enumerate}
\end{theorem}
\begin{proof}

Let $G=\Aut(L)$. The faces  of $\M$ which contain in-wedges form an orbit $W$ of $G$ in its action on faces. Faces not in $W$ must have pattern E. Let $f$ be one element of $W$. Whatever the pattern of $f$, that is the pattern of every face in $W$. Let $W_f$ be the set of in-wedges in $f$. 

We first note that $W_f$ is an orbit of some corner under the stabilizer $G_f$ in $G$, which is a subgroup of the dihedral group $D_p$.

If $f$ has pattern A, then every edge of $f$ is contained in two in-wedges in $f$ itself. Hence, none of its neighboring faces contain an in-wedge containing an edge of $f$. In particular it cannot be of type A and thus it must be of type E. But faces neighboring faces of type $E$, in turn, have in-wedges and so must be in $W$ which have pattern A. Thus, the map is face-bipartite, with one class of faces having pattern A, the other, pattern E.
  
If $f$ has pattern B, due the fact that every dart in contained in exactly one in-wedge, every neighboring face contains an in-wedge and must thus be of type B as well. By connectivity all faces are of type B. We may henceforth assume that no face is of type A or B.

Let $m$ be the length of a longest string $S$ of consecutive in-wedges in $f$. Note that $m$ cannot be greater than 2 if $f$ is not of type A; otherwise, there would be  an in-wedge adjacent to an out-wedge, and another in-wedge not adjacent to any out-wedge, contradicting that $L$ is transitive.

Now, if $m = 1$, then no two consecutive wedges of $f$ are in $L$. Since $f$ is neither of type $A$ or type $B$, there exists an edge of $f$ that belongs to no in-wedges of $f$, implying that the neighboring face, $f'$, of $f$ containing this edge, has two consecutive in-wedges. This contradicts the fact that $f$ and $f'$ are in the same orbit of faces. Therefore $m = 2$.

We then see there is an integer $k>1$ and integers $n_1, n_2, \ldots n_k$ such that the pattern of wedges around $f$ consists of a string $S_1$ of two in-wedges, followed by a string $T_1$ of $n_1$ out-wedges, followed by a string $S_2$ of two in-wedges, followed by a string $T_2$ of $n_2$ out-wedges, and so on. Since $L$ is transitive, there is a symmetry mapping the first in-wedge of $S_1$ to the first in-wedge of $S_2$. Such a symmetry cannot be a reflection, for it would map the last wedge in $S_1$ to the last wedge of $T_1$. Thus, it has to be a rotation and it must map each $S_i$ to $S_{i+1}$ and each $T_i$ to $T_{i+1}$. It follows that for some $n$, each $n_i$ is equal to $n$. Because half of all wedges in the map are in $L$, $n$ cannot be more than $m$.  So the remaining possibilities are

  %Suppose that a face $f$ has no consecutive in-wedges

%\micaelcomment{ Why is the following true? Couldn't we have $1$ black - $x$ whites - $1$ black - $y$ whites and then repeat this pattern?}
%If $f$ has some other pattern, then each face  in $W$ contains some $m$ consecutive corners in $L$,  followed by some $n$ consecutive corners not in $L$, then $m$ in $L$ and then $n$ not in $L$ and so on in a repeating pattern.  Because half of all wedges in the map are in $L$, $n$ cannot be more than $m$.  So the remaining possibilities are 

\begin{enumerate}
%\item $m = 1$, and so $n=1$, and in each face of $W$ in- and out-wedges alternate, giving pattern B.  As half of the wedges in $f$ are in $L$, and half of all wedges are in $L$, $W$ must consist of all faces.
\item $m = 2$ and $n = 1$, giving pattern C for face $f$.  Consider an edge  $e$ of $f$ belonging to two in-wedges.  In the face $f'$ across $e$ from $f$, the wedges containing $e$ must be out-wedges. Thus $f'$ is not in $W$, and so has pattern E.
\item  $m=2$ and $n = 2$.   Then in each face of $W$, pairs of in-wedges alternate with pairs of out-wedges: pattern D.  As before, the count shows that $W$ must be the set of all faces.
\end{enumerate}
This concludes the proof.
\end{proof}

\begin{remark}
Some information about the circuit structure $\CC(L)$ and local action of $\Aut(L)$ can be deduced in each of the four possible cases in the statement of Theorem \ref{th:faces}

\begin{enumerate}

\item If $L$ is of type (1), then $\M$ is face-bipartite and the circuits in the circuit structure $\CC(L)$ are the faces of one color in $\M$.  The local action is half-dihedral or half-cyclic. Further,  $L$ forms a corneration of type (2) in $\Pe(\M)$.

\item If $L$ is of type (2), then $G$ must be transitive on faces.   The elements of $\CC(L)$  are Petrie paths.  The local action is half-dihedral or half-cyclic. Further $L$ forms a corneration of type (1) in $\Pe(\M)$.

\item If $L$ is of type (3), then   $q$ must be divisible by 4, and  $G$ must not be transitive on faces.   The orbit $W$ contains at least $\frac{3}{4}$ of the faces.  The local action of $G$ is quarter dihedral.   Here, $L$ forms a corneration of type (3) in $\Pe(\M)$, as well.

\item If $L$ is of type (4), then $G$ must be transitive on faces.  The local action of $G$ is  half cyclic, and then $L$ forms a corneration of type (4) in $\Pe(\M)$.

\item In both (3) and (4), each circuit in $\CC(L)$  takes (locally) two right turns then two left, then two right, and so on.

\end{enumerate}
\end{remark}

%\begin{remark}
%If $L$ is of type (1) from the statement of the theorem, then the circuits in the circuit structure $\CC(L)$ are the faces of one color in $\M$.  Further,  $L$ forms a corneration of type (2) in $\Pe(\M)$.
%
%If $L$ is of type (2), then $G$ must be transitive on faces.   The elements of $\CC(L)$  are Petrie paths.    Further $L$ forms a corneration of type (1) in $\Pe(\M)$.
%
%
%If $L$ is of type (3), then   $q$ must be divisible by 4, and  $G$ must not be transitive on faces.   The orbit $W$ contains at least $\frac{3}{4}$ of the faces.  The local action of $G$ is quarter dihedral.   Here, $L$ forms a corneration of type (3) in $\Pe(\M)$, as well.
%
%If $L$ is of type (4), then $G$ must be transitive on faces.  The local action of $G$ is  half cyclic, and then $L$ forms a corneration of type (4) in $\Pe(\M)$.
%
%In both (3) and (4), each circuit in $\CC(L)$  takes (locally) two right turns then two left, then two right, and so on.
%\end{remark}

\begin{remark}
\label{rem:jvswedge}
The previous sections have given us an understanding of transitive 
$j$-cornerations  for which $j$ and $q$ are relatively prime, as these 
are equivalent to wedge-cornerations in $H_j(\M)$. From 
Remark \ref{rem:disconnect}, we see that even if $\gcd(q,j)>1$,  a  $j$-corneration 
$L$ of $\M$ is a  $1$-corneration of $H_j(\M)$; it must then exhibit one of the 
four patterns from Theorem \ref{th:faces}.
\end{remark}

\section{Symmetry-type graphs of cornerations}
\label{sec:STG}
%%**********************

%\primozcomment{We need to stick in the definition of reflexible, face reflexible and half reflexible maps and groups, and then observe thatthese correspond to certain symmetry types}
\color{Black}
 The topic of this section is the idea of a {\em symmetry-type graph}.  Informally, a symmetry-type graph of a group $G$ acting on a map $\M$ is a quotient of the flag graph of $\M$ with respect to $G$.   The symmetry-type graph  was first discussed in \cite{Dress}, under the name {\em Delaney symbol}.
A definition closer to the one we wish to present is found in \cite{HFOP}, where more (and earlier) references to this notion may be found.  In order to present  this  more rigourous definition, we need to introduce some preliminary ideas.

\subsection{Pre-graphs and diagrams}

We begin with the idea of a {\em pregraph}.   \color{Black}
  A pregraph is a quadruple $(V,D,\beg,\inv)$ where $V$ and $D$ are two disjoint non-empty sets and where $\beg \colon D \to V$ and $\inv  \colon D \to D$ are two mappings such that $\inv (\inv(x)) = x$ for every $x\in D$. 
 Elements of $V$ and $D$ are then called vertices and darts, respectively.
 An unordered pair $\{x,\inv(x)\}$, $x\in D$, is then called an edge of the pregraph.  The set of edges is denoted $E$.  If $x =\inv(x)$, the the edge $\{x,\inv(x)\}$ is said to be a {\em semiedge}. Similarly, if $\beg(x) = \beg(\inv(x))$ but $x\not=\inv(x)$, then the edge $\{x,\inv(x)\}$ is a {\em loop}.  All other edges are called {\em links}. A pregraph can be depicted by representing each vertex as a node (which might be a dot or a circle or a box). A link $\{x, \inv(x)\}$ is then
 represented as a line segment between the nodes representing $\beg(x)$ and $\beg(\inv(x))$, a loop $\{x,\inv(x)\}$ is drawn as a closed curve joining $\beg(x)$ to itself, while a semiedge $\{x\}$ is represented as a short line-segment emanating from $\beg(x)$.
An isomorphism from a pregraph $(V,D,\beg,\inv)$ to a pregraph $(V',D',\beg',\inv')$ is a bijection $\varphi \colon V\cup D \to V' \cup D'$, $\varphi(V) = V'$ and $\varphi(D) = D'$ such that $\varphi(\beg(x)) = \beg'(\varphi(x))$ and $\varphi(\inv(x))=\inv'(\varphi(x))$ for every $x\in D$.

Note that the notion of a pregraph generalises that of a graph. In particular, given a graph $\Gamma=(V,E,\partial)$, one can let
$D=\{\{v,e\} : e\in E, v\in\partial(e)\}$ be the set of darts of $\Gamma$, let $\beg(\{v,e\}) = v$ and let $\inv(\{v,e\}) = \{u,e\}$ where $\partial(e) = \{v,u\}$.
Then $(V,D,\beg,\inv)$ is a pregraph with additional property that it has no loops or semiedges. Conversely, one can easily see that every pregraph without loops or semiedges can be obtained from a graph in such a way. 

A mapping $c\colon E\to C$ from the edge-set $E$ of a pregraph to a set $C$ of things called colors is called an edge-coloring of a pregraph.
Such an edge-coloring is {\em locally bijective} provided that the restriction of $c$ to the set of edges incident to a vertex $v$ is bijective onto $C$
for every $v\in V$. In particular, existence of a locally bijective edge-coloring implies that the valence of every vertex is equal to the cardinality of $|C|$.
A pair $(\Gamma,c)$ where $\Gamma$ is a pregraph and $c\colon E\to C$ is a locally-bijective edge-coloring of $\Gamma$ is called an edge-colored pregraph (or, to be more precise, an edge-$|C|$-colored pregraph).
%If the image of an edge-coloring is of cardinality $k$, then we will call $c$ a $k$-coloring.
 If $(\Gamma,c)$ and 
$(\Gamma',c')$ are two edge-colored pregraphs and $\varphi \colon \Gamma \to \Gamma'$ is a pregraph isomorphism such that
$c(e) = c'(\varphi(e))$ for every edge $e$ of $\Gamma$, the we call $\varphi$ an isomorphism of edge-colored graphs.

Finally, for the purposes of this paper, we define a {\em vertex-shaped edge-$3$-colored pregraph} (or simply a {\em diagram}, for short) to be
a edge-$3$-colored pregraph together with a labelling $V \to \{{\rm B}, {\rm O}\}$, where vertices assigned the label $B$ will be called {\em boxes}
and depicted as squares, while those assigned the label O will be called {\em ovals} and will be depicted as such.
Two diagrams are then isomorphic if they are isomorphic as edge-colored pregraphs via a vertex-label preserving isomorphism.

\subsection{Symmetry-type graphs}
The {\em flag (pre)graph} of a map $\M$ with the flag-set $\Omega$ is a edge-colored pregraph $(\Gamma,c)$,
where the vertex-set of $\Gamma$ is $\Omega$, and a dart (of color $i$) is an ordered pair $x=(f,r_if)$ for some $i\in\{0,1,2\}$,
where $\beg(x) = f$ and $\inv(x) = (r_if,f)$. Let $\Delta_i$ denote the set of all darts of color $i$ and let $\Delta = \Delta_0\cup\Delta_1\cup\Delta_2$.
Note that the colors of $x$ and $\inv(x)$ are the same so we may define the color $c(e)$ of  the edge $e=\{x,\inv(x)\}$ to be their common color.
An important observation is that the edges of colors $0$ and $2$ always induce closed walks of length dividing $4$.

\begin{figure}[hhh]
\begin{center}
\epsfig{file=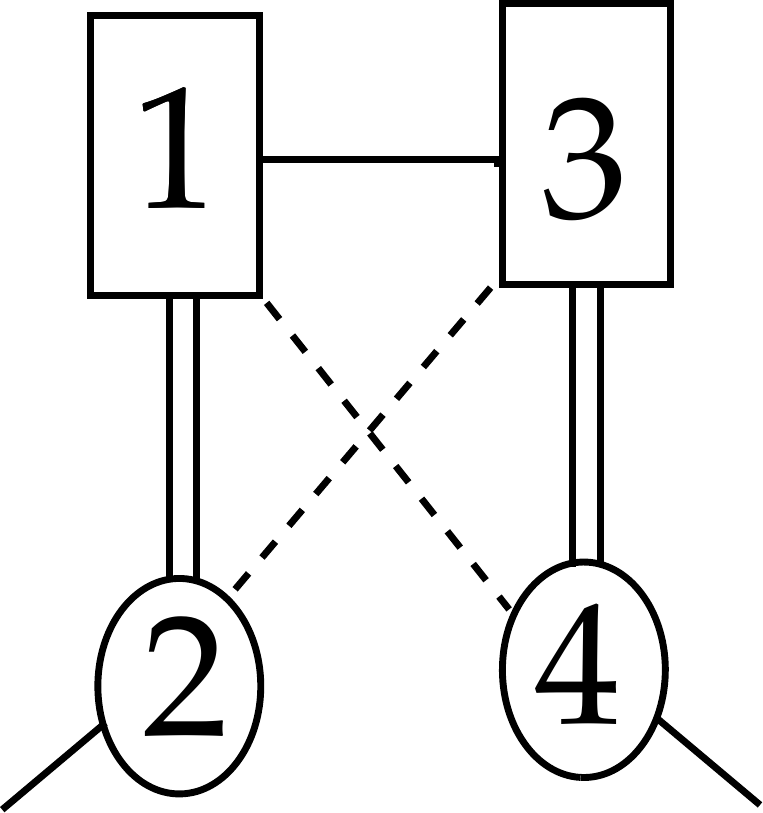,height=30mm}
\caption{A symmetry-type graph (in fact, a diagram) of a map with respect to a group of symmetries acting on a $1$-corneration. The $0$-edges are depicted with dotted lines, $1$-edges with single solid lines and $2$-edges with double solid lines. Boxed vertices $1$ and $3$ correspond to orbits of in-wedges while the ovals $2$ and $4$ correspond to orbits of out-wedges.}
\label{fig:STG8}
\end{center}
\end{figure}

Assuming the above notation, we can now define
the {\em symmetry-type graph} $\M/G$ of the map $\M$ with respect to a group of symmetries $G$ as a pregraph
$(\Omega/G,\Delta/G,\beg',\inv')$, where $\Omega/G$ and $\Delta/G$ are the orbits of the natural action of $G$ on $\Omega$ and $\Delta$,
respectively, and where $\beg'(x^G) = \beg(x)^G$ and $\inv'(x^G) = \inv(x)^G$. That is, the symmetry type graph is the quotient (as defined in \cite{GenVolt}, for example) of the flag-graph of $\M$ with respect to the action of $G$.
Finally, since the action of $G$ on the flag-graph clearly preserves the colors of the edges, this allows us to define a color of an edge $e^G=\{x^G,\inv'(x^G)\}$ to be that of the edge $e$ in the flag-graph. Note that this coloring is locally-bijective and hence the symmetry type graph is an
edge-$3$-colored pregraph.
The symmetry type graph can be visualised by drawing as a pregraph with $0$-, $1$- and $2$-edges depicted as dotted, single and double lines,
as in Figure~\ref{fig:STG8}. 

Furthermore, when given a $1$-corneration $L$ of a map $\M$, we may assign the label B (box) to every flag of $\M$ contained in an in-wedge of $L$,
and the label O (oval) to every flag contained in an out-wedge. In this way, a symmetry type graph $\M/G$ with respect to some subgroup $G\le\Aut(L)$
is a diagram. Consider the sample diagram in Figure~\ref{fig:STG8}.

Note how much about $\M$ we can deduce from $\M/G$:
\begin{itemize}
\item The consecutive flags of one face are 0-adjacent and 1-adjacent.  Because the 0-edges and 1-edges of the symmetry-type graph join all four nodes in this diagram, every face contains a flag of each orbit.  Thus the group acts transitively on faces.\color{Black}
\item Looking at those edges in slightly more detail,  we see that a sequence of 0-1-paths returns to the starting flag in four steps, we deduce that for each face, the group contains a  rotation by four steps about that face.
\item  Because of the semi-edge at node 2, each flag of orbit 2 is 1-adjacent to another flag of the same orbit.  \color{Black}
\item As the flags around a vertex are joined by 1- and 2- adjacencies, and the 1- and 2- edges in the diagram connect all nodes, the group acts transitively on vertices as well. Thus $G_v$ is $QD$.
\item Two flags in a wedge are 1-adjacent.  If they are in an in-wedge (rectangular nodes), one is in orbit 1, the other (if there is another) in orbit 3.  Thus the group acts transitively on in-wedges.
\item Of the two flags   in an out-wedge (oval nodes), if one is in orbit 2, the other in orbit 2, as well. If one is in 4, so is the other. Thus the group has two orbits on out-wedges.
It follows that the group must have a reflection at each such corner, and similarly, a reflection at each flag of orbit 4.
Therefore faces have pattern B.
\end{itemize}

\begin{remark}\label{rem:diag}
If $\M$ is a map with a $G$-transitive wedge-corneration, then the following holds for the symmetry type graph $\M / G$:
\begin{enumerate}
\item $\M /G$ is connected and has $2$ or $4$ nodes;
\item an edge of color $2$ must connect nodes of different  shapes (and thus there no semi-edges of color $2$);
\item an edge of color $1$ must connect nodes of the same type (and thus there are no full edges of color $1$ if $\M/G$ has only two vertices, but only semi-edges);
\item two distinct rectangular nodes must be connected by an edge of color $1$;
\item a walk of length $4$ along the edges and semi-edges of $\M / G$ that alternates colors $0$ and $2$ must be closed;

\end{enumerate}
\end{remark}

\begin{figure}[hhh]
\begin{subfigure}{1.4in}
\begin{center}
\includegraphics[height=19mm]{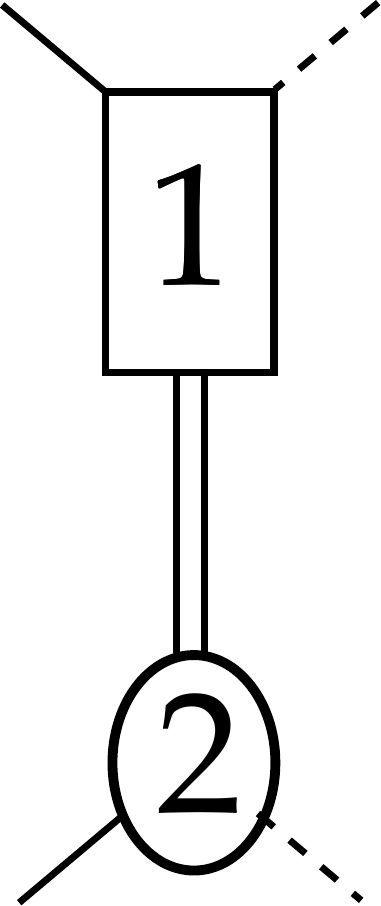}
\caption{}\label{fig:X1}
\end{center}
\end{subfigure}
~
\begin{subfigure}{1.4in}
\begin{center}
\includegraphics[height=17mm]{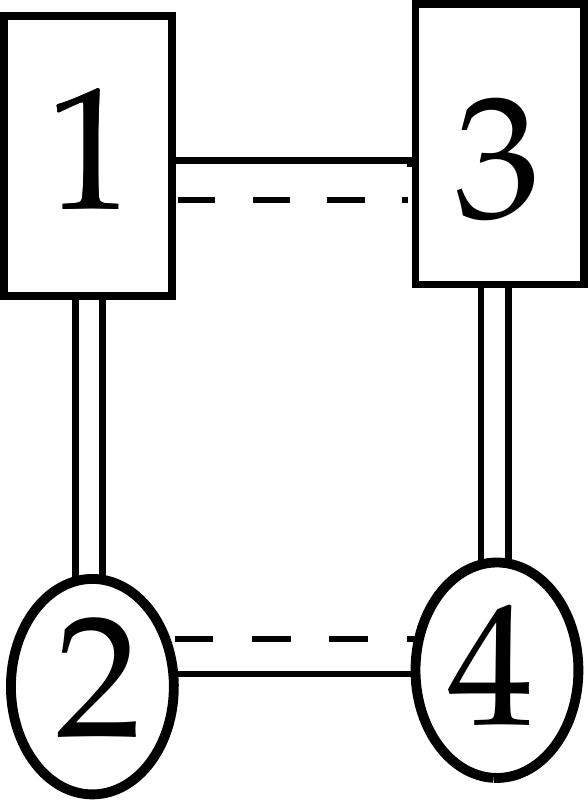}
\caption{}\label{fig:X2}
\end{center}
\end{subfigure}
 ~
\begin{subfigure}{1.4in}
\begin{center}
\includegraphics[height=19mm]{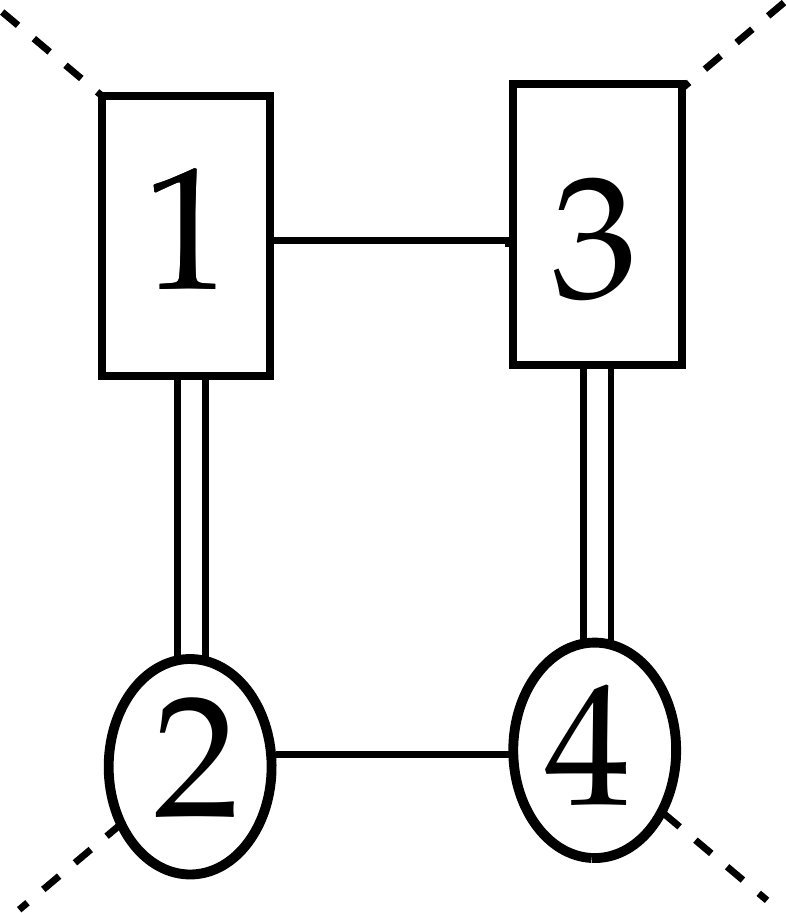}
\caption{}\label{fig:X3}
\end{center}
\end{subfigure}
 ~
\begin{subfigure}{1.4in}
\begin{center}
\includegraphics[height=18mm]{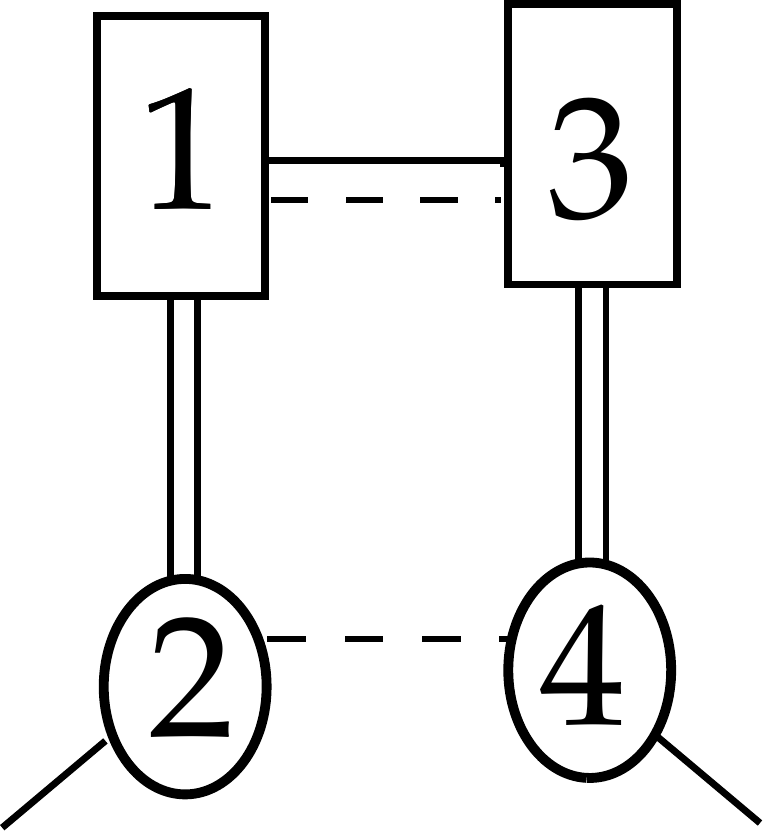}
\caption{}\label{fig:X4}
\end{center}
\end{subfigure}

\begin{subfigure}{1.4in}
\begin{center}
\includegraphics[height=19mm]{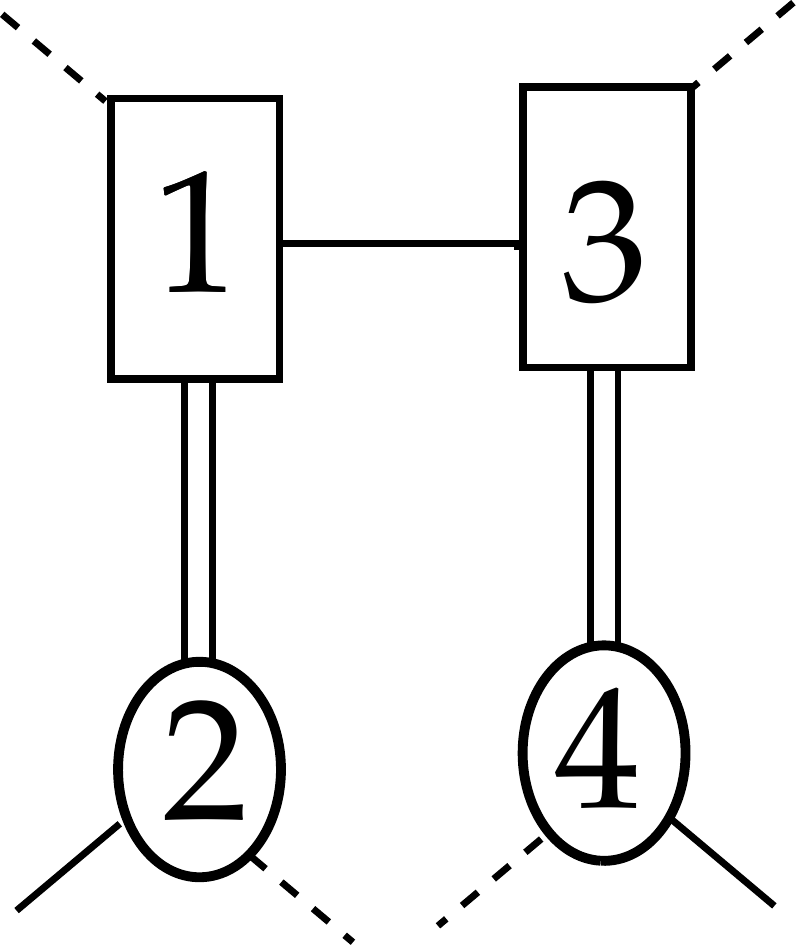}
\caption{}\label{fig:X5}
\end{center}
\end{subfigure}
~
\begin{subfigure}{1.4in}
\begin{center}
\includegraphics[height=19mm]{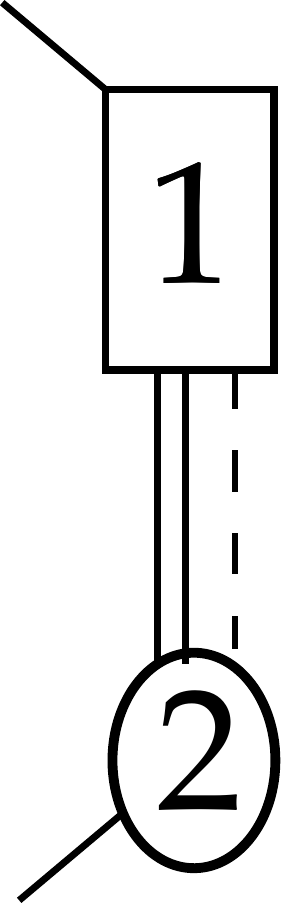}
\caption{}\label{fig:X6}
\end{center}
\end{subfigure}
 ~
\begin{subfigure}{1.4in}
\begin{center}
\includegraphics[height=17mm]{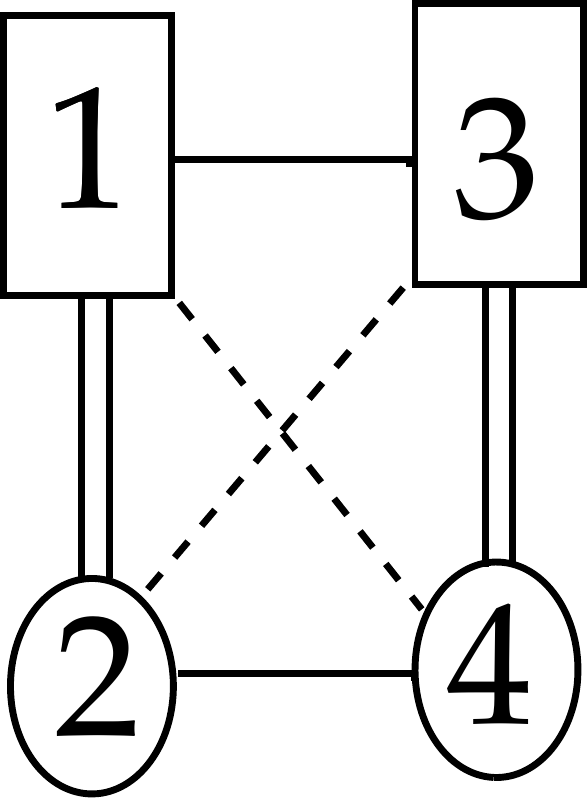}
\caption{}\label{fig:X7}
\end{center}
\end{subfigure}
 ~
\begin{subfigure}{1.4in}
\begin{center}
\includegraphics[height=19mm]{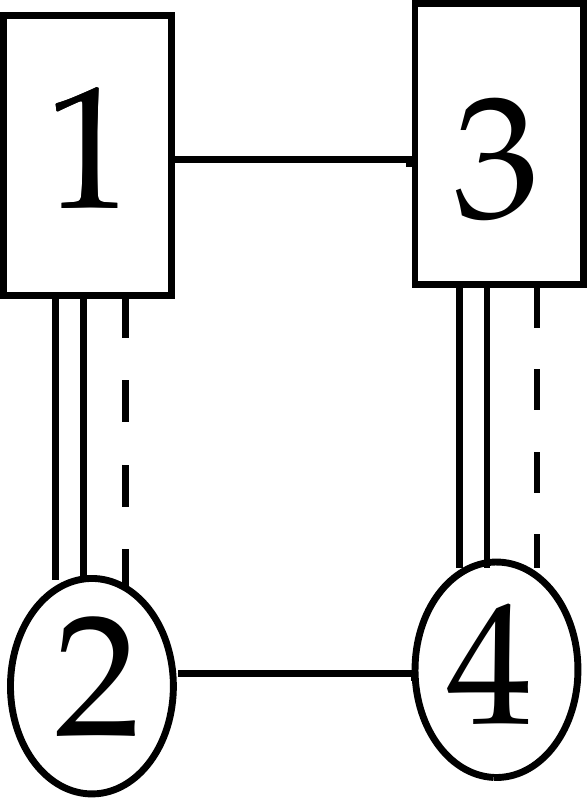}
\caption{}\label{fig:X8}
\end{center}
\end{subfigure}

\begin{subfigure}{1.4in}
\begin{center}
\includegraphics[height=17mm]{STG8.pdf}
\caption{}\label{fig:X9}
\end{center}
\end{subfigure}
~
\begin{subfigure}{1.4in}
\begin{center}
\includegraphics[height=18mm]{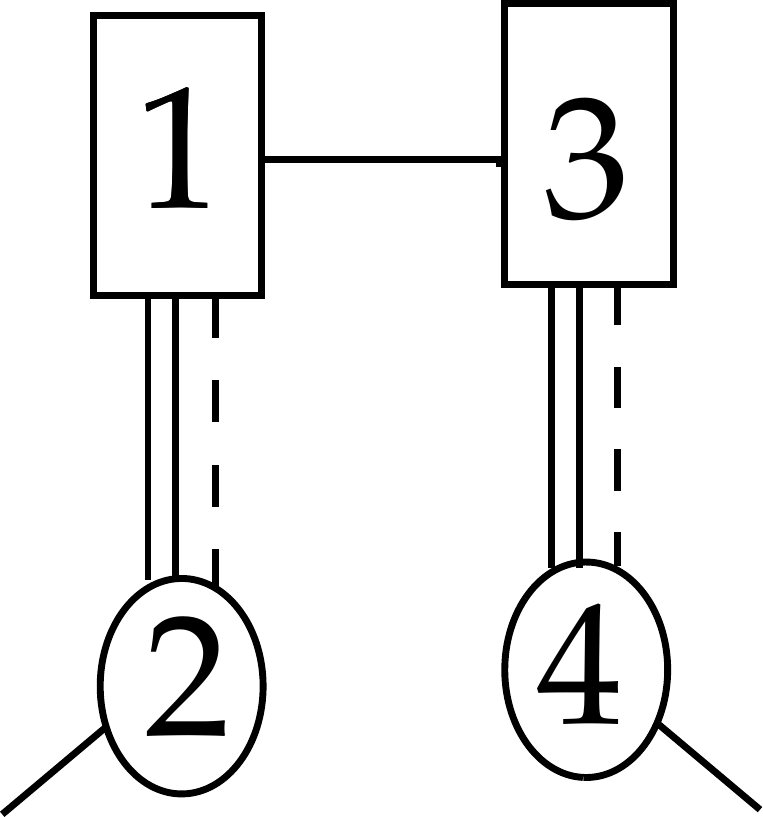}
\caption{}\label{fig:X10}
\end{center}
\end{subfigure}
 ~
\begin{subfigure}{1.4in}
\begin{center}
\includegraphics[height=19mm]{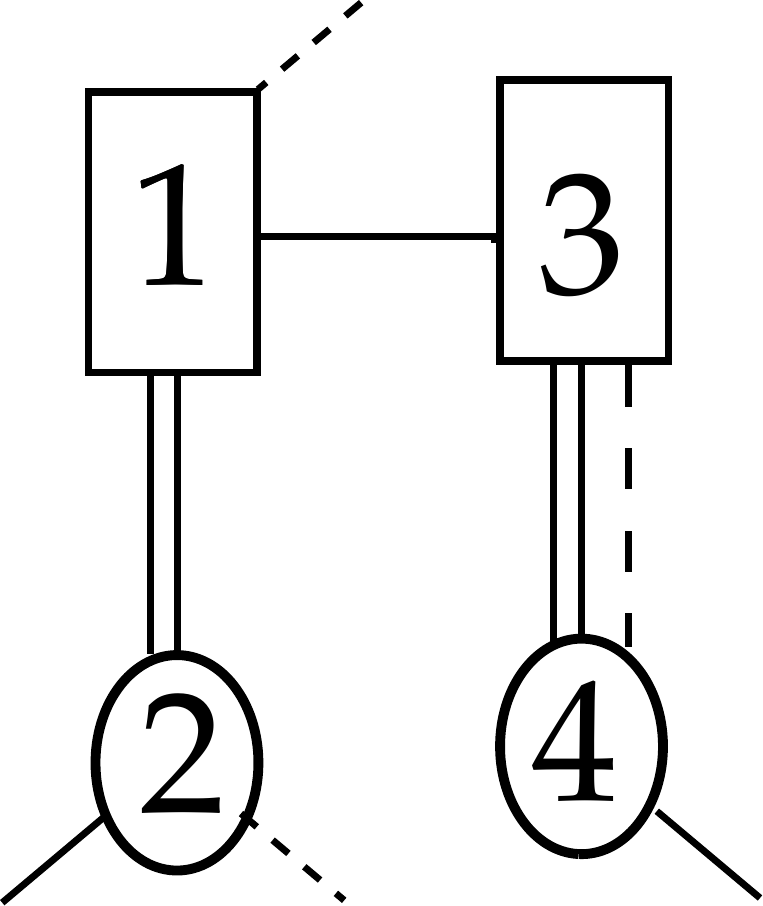}
\caption{}\label{fig:X11}
\end{center}
\end{subfigure}
 ~
\begin{subfigure}{1.4in}
\begin{center}
\includegraphics[height=19mm]{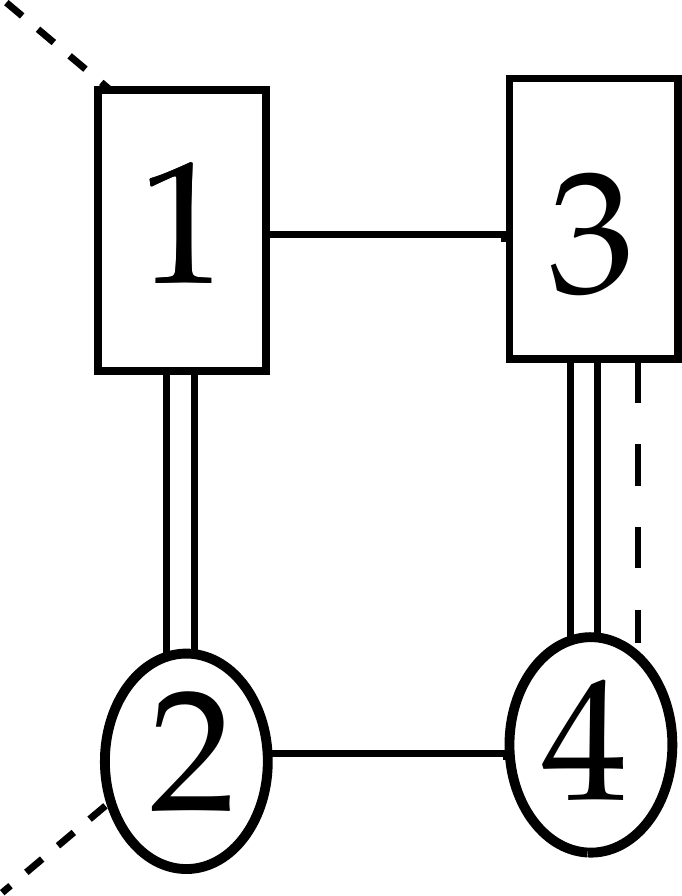}
\caption{}\label{fig:X12}
\end{center}
\end{subfigure}

\caption{All symmetry-type graphs of transitive cornerations}
\label{fig:STGall}
\end{figure}

There are 12 possible symmetry-type graphs for a transitive corneration, and these are shown in Figure \ref{fig:STGall}. From these diagrams one can easily deduce information on the number of orbits on vertices, edges and faces of the corresponding quotient group, as well as the type of $1$-corneration (see Theorem \ref{th:faces}) and the local action. This is summarized in table \ref{tab:CornGps}

\begin{table}[hhh]
\begin{center}
\begin{tabular}{|c|c|c|c|c|c|}
\hline
STG  &  V orbs &E orbs  &  F orbs&Pattern&local type\\
\hline
 
(a) &1&1&2&A, E&HD\\
\hline
(b) &1&1&1&A, E&HC\\
\hline
(c) &1&2&2&A, E&HC\\
\hline
(d) &1&1&2&A, E&QD\\
\hline
(e) &1&2&3&A, E&QD\\
\hline

(f) &1&1&1&B&HD\\
\hline
(g) &1&1&1&B&HC\\
\hline
(h) &1&2&2&B&HC\\
\hline
(i) &1&1&1&B&QD\\
\hline
(j) &1&2&2&B&QD\\
\hline

(k) &1&2&2&C, E&QD\\
\hline
(l) &1&2&1&D&HC\\
\hline

\end{tabular}
\caption{Transitive corneration groups. Each line represents one of the 12 possible symmetry type graphs of transitive wedge cornerations, with
the letters (a)--(l) in column STG corresponding to the 12 diagrams in
 Figure~\ref{fig:STGall}. The next three columns (V orbs, E orbs and F orbs) show the number of orbits on vertices, edges and faces, respectively of the transitive corneration group yielding the diagram. The column `Pattern' shows the face patterns (as given in Figure~\ref{fig:XABCD}) appearing in the corneration. The last column,  labelled local type, shows the type of the local action as determined in Theorem~\ref{th:locgps}.}
\label{tab:CornGps}
\label{Tab2}
\end{center}
\end{table}

\begin{theorem}
\label{the:class}
 A diagram is a symmetry type graph of a map $\M$ with respect to a group $G$ acting transitively on a $1$-corneration of $\M$
 if and only if it is isomorphic to one of the diagrams in Figure \ref{fig:STGall} and this happens if and only if it satisfies the  requirements on $\M/G$ in Remark \ref{rem:diag}.
%A diagram appears in Figure \ref{fig:STGall} exactly when it satisfies the  requirements on $\M/G$ in Remark \ref{rem:diag} and this happens exactly when it is  the symmetry-type graph  of a transitive corneration  of some map.
\end{theorem}

\begin{figure}[hhh]
\begin{center}
\epsfig{file=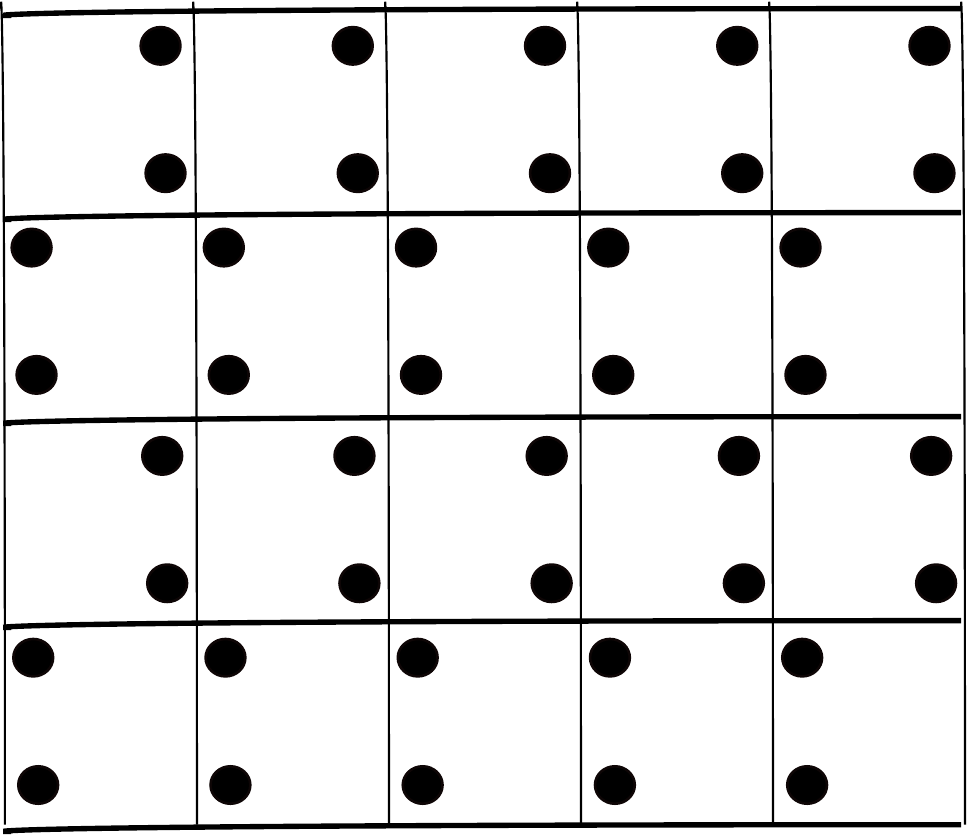,height=30mm}
\caption{A corneration in a torus map}
\label{fig:4445}
\end{center}
\end{figure}

\begin{proof}

The reader may convince themselves that every diagram in Figure \ref{fig:STGall} satisfies the requirements of Remark \ref{rem:diag} and that every such diagram appears in  Figure \ref{fig:STGall}.  We then need only to show that each of these diagrams actually occurs as the symmetry-type graph of a transitive corneration. In what follows, we will exhibit a map $\M$ admitting ten different transitive $1$-cornerations, each realizing one of the diagrams (a)--(j) of Figure \ref{fig:STGall}. This map $\M$ will be defined in term of the opposite operator, $\textrm{opp}$.  Recall that if $\N$ is a map, then $\textrm{opp}(\N)$ is defined as the dual of the Petrie of the dual of $\N$. 

As was explained in \cite{Wop}, $\textrm{opp}(\N)$ can also be obtained by separating the faces of $\N$ and glueing them again in such a way that if two faces $F_1$ and $F_2$ shared an edge with endvertices $u$ and $v$ then in $\textrm{opp}(\N)$ the segment from $u$ to $v$ in $F_1$ gets identified with the segment from $v$ to $u$ in $F_2$. In that sense, we can say that $\N$ and $\textrm{opp}(\N)$ have the same set of faces, and two faces are adjacent in $\N$ if and only if they are adjacent in $\text{opp}(\N)$. In particular if $\N$ is face-bipartite then so is $\text{opp}(\N)$.

Consider the graph with vertex-set $\ZZ_4 \times \ZZ_4$ with `horizontal' edges of the form $(x,i) \sim (x,i+1)$ and `vertical' edges of the form $(x,i) \sim (x+1,i)$. By glueing a disc along each $4$-cycle of alternating horizontal and vertical edges, one obtains a face-bipartite $4$-valent map $\N$  on the torus with faces of length four and Petrie paths of length $8$/ This map is sometimes denoted as $\{4,4\}_{4,0}$ (see \cite{CoxMos}) and is named $\textrm{RM}[32,3]$ in the census of rotary maps \cite{mapscensus}.  

Now consider the opposite map $\M = \textrm{opp}(\N)$  (this map is named $\textrm{RM}[32,4]$ in \cite{mapscensus} and it sits on an orientable surface of genus $5$). %The automorphisms groups of $\N$ and $\M$ are isomorphic and act transitively on the flags of the two maps.
It is an $8$-valent face-bipartite map with Petrie paths of length $4$.  If we think of the faces as being colored red and green, and if we let $L$ be the set of all wedges of green faces, this gives us a corneration of pattern $A$ (in the green faces). Since $\Aut(\M)$ acts transitively on the flags of $\M$, the group $G := \Aut(L)$ has index $2$ in $\Aut(\M)$ and acts transitively on flags in green faces as well as on flags in red faces, implying that $L$ is a transitive corneration. In particular, $G$ has two orbits on wedges and its symmetry-type graph is diagram (a) in Figure \ref{fig:STGall}. Further, one can readily check with a computer package such as GAP \cite{gap} or Magma \cite{magma} (and using the census \cite{mapscensus}) that $G$ has four subgroups, each of index 2, each having 4 orbits of flags, and each are transitive on the wedges in $L$.  Their symmetry-type graphs are diagrams (b)--(e) in  Figure \ref{fig:STGall}.
\color{black}

The map $\M$ is isomorphic to $\Pe(\M)$.  The corneration $\Pe(L)$ has pattern $B$ and the operator $\Pe$ turns diagram (a) into (f), (b) into (g) and so on.   Thus diagrams (a)--(j) have realizations on $\M$.

Diagram (k) has many realizations, but the easiest to see are in the geometric anti-prisms.   As a map on the sphere, the $n$-antiprism is formed by two $n$-gonal faces joined by a band of $2n$ triangles.  The corneration $L_n$ contain none of the wedges in the $n$-gons, but does contain the two wedges of each triangle that are incident with the edge joining the triangle to an $n$-gon.  The case $n = 4$ is shown in Figure \ref{fig:AntiP}.

\begin{figure}[hhh]
\begin{center}
\epsfig{file=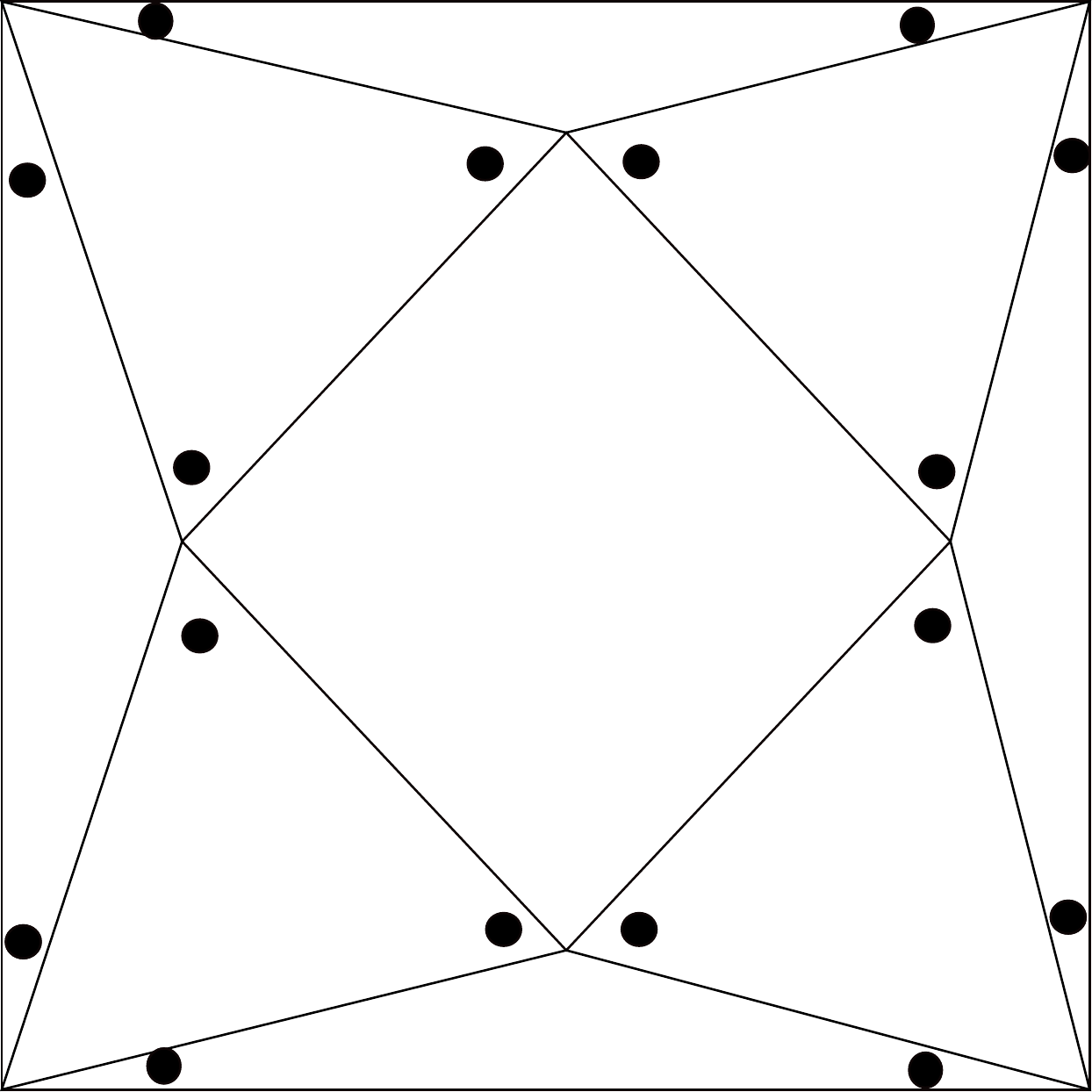,height=30mm}
\caption{A corneration in the $4$-antiprism}
\label{fig:AntiP}
\end{center}
\end{figure}

Diagram (l), also, has many realizations, and the easiest to see in this case are on the torus.   Let the map $\M$ on the torus be formed from a $2m\times n$ rectangle of squares, with opposite sides identified in the usual way.  The corneration $L$ contains wedges incident with one of the two vertical edges in each square.  The choice of which of those two edges to use is made consistently along each row of the grid but alternatingly along each column.  The special case with a $4\times 5$ rectangle is shown in Figure \ref{fig:4445}. 
\end{proof}
%% or in Figure \ref{fig:44452}

%\begin{figure}[hhh]
%\begin{center}
%\epsfig{file=44452.pdf,height=25mm}
%\caption{A corneration in a torus map}
%\label{fig:44452}
%\end{center}
%\end{figure}

\color{Black}
\section{Symmetric Cornerations}
\label{sec:SymCor}
For many applications, we are interested in the case in which $\Aut(L)$ acts transitively on the darts of the maps $\M$; these are the symmetric cornerations.  We want to determine the maps in which these occur and how to find $L$ is such a map. Recall that a dart can be seen as an orbit of flags under $r_2$. Therefore, if $\Aut(L)$ is arc-transitive, then the subgraph of the corresponding symmetry-type graph induced by edges of type $2$ must have a single connected component (or, if you will, the edges of type $2$ are a spanning tree in $\M/G$). Scanning Figure \ref{fig:STGall}, we see that only Diagrams (a) and (f) satisfy this condition, and thus any symmetric corneration must have one of those two as its symmetry-type graph.

	Notice that in Diagram (a), all the flags of any one face (that is, an orbit under $r_0$ and $r_1$) are in one orbit, and therefore, the group must contain all reflections fixing that face.  These are all the reflections whose axes pass through the center of that face.

%%\primozcomment{We need to prove the following:}

%\begin{theorem}
% Let $L$ be a $j$-corneration of a map $\M$. Then $L$ is symmetric if and only if $\M$ or its Petrie has a subgroup that acts half reflexibly on $\M$.
%\end{theorem}

\begin{theorem}
\label{the:symodd}
Let $\M$ be a map of valence $q$ and let $j < \frac{q}{2}$. Then, $\M$ admits a symmetric $j$-corneration if and only if $j$ is odd and $\Aut(\M)$ has a subgroup that
 acts half-reflexibly on $\M$ or $\Pe(\M)$.
\end{theorem}

\begin{proof}
First let $L$ be a symmetric $j$-corneration of a map $\M$.  Recall that this means that $\Aut(L)$ is transitive on the darts of $\M$, and that, by Lemma \ref{lem:DTsym}, $j$ is odd.  Then the local corneration is the standard odd local corneration, 
and so the
stabilizer of a vertex must contain the reflections we called $\mu_i$ in section \ref{sec:local} exactly when $i$ is odd.  These are reflections about each axis joining a vertex to a face-center.  

For an edge 
$e=\{u, v\}$, the corners $L(e,u)$ and $L(e,v)$ either 
align convexly or with inflection, consistently throughout the map. Without loss of generality, suppose $\M$ (and not its Petrie) is the map for which those alignments are all convex.  Then the element of $\Aut(L)$ which reverses 
an edge in $\M$ must be a reflection about the axis joining a face-center to an edge-midpoint.  Thus $\Aut(L)$ acts half-reflexibly on $\M$, proving the first assertion of the theorem.

	Conversely, suppose that $\M$ is a map, that $G \le \Aut(\M)$ is a group acting half-reflexibly on $\M$, and that $j$ is odd.    Then $G$ has two orbits on flags, and we may think of those orbits as colors, red and green. Then for any one face, its flags are all of one color, while any two adjacent faces have opposite colors.

 Wedges are orbits under $r_1$, and so each wedge is a single color.  Because $j$ is odd, the interior boundary wedges of any    $j$-corner are the same color.   Define $L_R$ to be the set of all $j$-corners of $\M$ whose interior boundary  wedges are red and similarly define $L_G$ to be the set of all $j$-corners of $\M$ whose interior  boundary wedges are green.   Then $L_R$ and $L_G$ are two cornerations each preserved by $G$, and thus are symmetric.
\end{proof}

\begin{corollary}
Let $\M$ be a reflexible map of valence $q$ and $j<\frac{q}{2}$ a positive integer. Then $\M$ admits a symmetric $j$-corneration if and only if $j$ is odd and at least one of $\M$ and
$\Pe(\M)$ is face-bipartite.
\end{corollary}

\begin{proof}
This follows directly from Theorem~\ref{the:symodd} and Lemma~\ref{lem:HR}.
\end{proof}

%%\primozcomment{Can we say something about how many symmetric $j$-cornerations we get under the above circumstances?}

\color{black}

\section{Split graphs}
\label{sec:SG}

In this section, we first define $\SSS(L,K)$ for a corneration $L$ and a set $K$ of corners.   This is a the generalization of the Split Graph of a cycle decomposition.  We give examples $\A(L), \B(L), C_i(L), C_x(L)$ of the $\SSS$ construction and indicate when these are vertex-transitive.   We then discuss (with proof)  which of the graphs are locally connected (and thus connected),  and the valence of each.

%\subsection{Definitions}\label{sec:SGdef}
\begin{definition}\label{def:SS}
 Let $L$ be a corneration of $\Gamma$, and let $K$ be any collection of corners in $\Gamma$ disjoint from $L$.  
 Let the split graph $\SSS(L, K)$  be the simple graph whose vertices are the elements of $L$ and whose edges are of the two following kinds:
\begin{enumerate}
\item {\em `Old'} edges  of $\SSS(L, K)$ are all $\{c_1, c_2\}$, where $c_1$ and $c_2$ are corners in $L$ such that $c_1\cap c_2 = \{e\}$ for some edge $e$ of $\Gamma$.
\item {\em `New'} edges of $\SSS(L, K)$ are all pairs of the form $\{L(\{e_1,v\}), L(\{e_2,v\})\}$, where $\{e_1, v, e_2\} \in K$.
\end{enumerate}
\end{definition}

\begin{proposition}\label{pr:trans}
If $L$ is a transitive corneration of a map $\M$ and 
 a set  $K$ of corners disjoint from $L$ is invariant under a subgroup $G$ of $\Aut(L)$ acting transitively on $L$,
 then $\SSS(L, K)$ is a vertex-transitive graph.
\end{proposition}

\begin{proof}
\color{Black}
It should be clear that $H=\Aut(L)$ is a subgroup of $G = \Aut(\Gamma)$, where the latter is considered a group of permutations on the vertices and edges of $\Gamma$.  Because $H$ acts  on $L$, it acts  on the vertices and old edges of $\SSS(L, K)$.   Because $H$ acts on $K$, it acts on new edges of $\SSS(L, K)$ as well.

Because the action on $L$ is transitive, the action of $H$ on  $\SSS(L, K)$ is vertex-transitive, as required.
\end{proof}

\color{Black}
Let $L$ be a corneration in a graph $\Gamma$ and let $K$ be a set of corners in $\Gamma$ disjoint from $L$.
Then the split graph $\SSS(L, K)$ is called {\em locally connected} provided that for every $v\in \V(\Gamma)$
the subgraph of $\SSS(L,K)$ induced by all the corners in $L$ containing $v$ is connected.
Note that if the split graph $\SSS(L,K)$ is locally connected, then it is connected, but the converse does not hold:  it can happen that $\SSS(L,K)$ is connected even though it is not locally connected.
The smallest example appears to come from  3-uniform cornerations on a reflexible map of type $\{3,12\}_9$ on 162 edges and its Petrie.  This map is item $[162,56]$ in the on-line census of rotary maps on (the section on regular maps) \cite{mapscensus}  and its Petrie  is item  $[162,69]$ there.

\subsection{Four constructions}

We present here four examples of graphs constructed using Definition \ref{def:SS}, and we prove that each is vertex-transitive when $L$ is transitive.  Then in the following subsection,  we will determine  for each, its valence and  criteria  for it to be locally connected
.
 \begin{definition}
If $L$  is a $j$-uniform corneration of a map $\M$ for $j<\frac{q}{2}$, let  $\A(L)=\SSS(L, K)$, where $K$ is the $j$-complement of $L$ in $\M$.
 \end{definition}

\begin{definition}
If $L$  is a $j$-uniform corneration of a map $\M$, $2 \leq j \leq \frac{q}{2}$, and $K$ is the set of all wedges in $\M$, let $\B(L)=\SSS(L, K)$.
 \end{definition}

\begin{definition}
If $L$  is a $j$-uniform corneration of a map $\M$, $2 \leq j < \frac{q}{2}$, and $K_i$ is the set containing every wedge in $\M$ which is an interior boundary wedge of some corner of $L$, let $\C_i(L)=\SSS(L, K_i)$.
 \end{definition}

\begin{definition}
If $L$  is a $j$-uniform corneration of a map $\M$, $2 \leq j < \frac{q}{2}$, and $K_x$ is the set containing every wedge in $\M$ which is an exterior boundary wedge of some corner of $L$, let $\C_x(L)=\SSS(L, K_x)$.
 \end{definition}

\begin{theorem}
If $L$ is a transitive corneration, then each of the graphs $\A(L)$, $\B(L)$, $C_i(L)$, $C_x(L)$ is vertex-transitive.
\end{theorem}
\begin{proof}
In each case, the group $\Aut(L)$ preserves the corresponding set $K$, and so 
Proposition \ref{pr:trans} applies.
\end{proof}

\subsection{Local Connectivity and Valence of Split Graphs}   \color {Black}
Local connectivity is decided in the local corneration, and we use and enlarge on the notation from Section \ref{sec:local}:
\begin{itemize}
\item $L_v$ is the collection of edges and corners containing the vertex $v$.  
\item Edges are numbered  $0,1,2,\dots, q-2, q-1$ in order consecutively around $v$
\item A corner is just a pair  $\{a, b\}$
\item If the width of $L$ is $j$, $c_i$ stands for the corner $\{i, i+j\}$.
\end{itemize}
We first deal with the corneration of all straight corners.  Here, only the definition for $\B(L)$ applies.

\begin{theorem}\label{th:straight}
Let $\M$ be a dart-transitive map, and let $L$ be the set of straight corners.  Then $L$ is transitive and $\B(L)$ is vertex-transitive and locally connected.  Its valence is $4$ unless $q = 4$, in which case  the valence is $3$.

\end{theorem}

\begin{proof}
The straight corners at $v$ are  all $c_i = \{i, i+\frac{q}{2}\}$. The wedges sharing an edge  with the  straight corner $c_0$ are 
$\{0,1\}, \{\frac{q}{2}-1, \frac{q}{2}\},\{\frac{q}{2}, \frac{q}{2}+1\}, \{q-1,0\}$.  The first and third of these join $c_0$ to $c_1$ and  the second and fourth join $c_0$ to  to $c_{-1}$   Similarly, each $c_i$ is joined to $c_{i-1}$ and $c_{i+1}$, and so $\B(L) $ is locally connected.  It is of valence $4$ unless $c_1 = c_{-1}$, which happens if and only if $q = 4$
\end {proof}

%%\primozcomment{Note that $\A(L)$ is locally connected if and only if $j$ is relatively prime to the valence $q$.}

 As we have dealt with straight cornerations in Theorem \ref{th:straight}, we can now assume that the width $j$ of any transitive corneration is strictly less than $\frac{q}{2}$, so that all four of  $\A(L)$, $\B(L)$, $ C_i(L)$, $C_x(L)$ are defined.

 Of these, we now consider the graph $\A(L)$:

\begin{theorem}\label{th:Aval}
If $L$ is a transitive corneration of width $j$, then $\A(L)$ is locally connected if and only if $\gcd(q,j) = 1$.  Its valence is $4$, unless $j = \frac{q}{4}$, in which case the valence is 3.
\end{theorem}
\begin{proof}
We see from Remark \ref{rem:comp} that $\A(L)$ is locally connected exactly when $\gcd(q,j) = 1$.  
Assuming, as we may, that $c_0\in L$, and letting $K$ be the $j$-complement of $L$, we have that $c_0 \in L, c_j \in K, c_{2j}\in L, \dots $.
  Thus, corner $c_{0}$ is connected to $c_{2j}$ and  to $ c_{-2j}$.  
Then $c_{0}$ is connected to two elements of $L$ (and so has valence $4$ in $\A(L)$)  unless $2j = -2j$.  
This last can happen only if $4j = 0$, and then, the valence of $v$ is $3$.
\end{proof}

We will present results about the remaining graphs,  $\B(L), C_i(L), C_x(L)$ together, as the wedge-set in the construction of $\B(L)$ is the union of those for $C_i(L)$ and  $ C_x(L)$. but we will separate the presentation into two cases: the local corneration is the standard odd local corneration and the local corneration is the standard even  local corneration.

%%\primozcomment{Note that $\B(L)$ is always locally connected and thus connected.}
\begin{theorem}
\label{th:BCodd}
If $L$ is a transitive corneration of odd width $3\le j<\frac{q}{2}$, then:
\begin{enumerate}
\item $\B(L)$ is locally connected,
\item $\C_i(L)$ is locally connected exactly when $\gcd(q, j-1) = 2$,
\item $\C_x(L)$ is locally connected exactly when $\gcd(q, j+1) = 2$.
\end{enumerate}
Furthermore:
\begin{enumerate}
\item[(5)] $\B(L)$ has valence $6$ unless $q$ is divisible by $4$ and $j = \frac{q}{2}-1$, in which case it has valence $5$,
\item[(6)] $\C_i(L)$ and $C_x(L)$ each have valence $4$ unless $q$ is divisible by $4$ and $j = \frac{q}{2}-1$, in which case they have valence $3$.
\end{enumerate}
\end{theorem}
\color{Black}

%\begin{theorem}\label{th:BCodd}
%If $L$ is a transitive corneration of odd width $3\le j<\frac{q}{2}$, then $\B(L)$ is locally connected, $\C_i(L)$ is locally connected exactly when $\gcd(q, j-1) = 2$, while $\C_x(L)$ is locally connected exactly when $\gcd(q, j+1) = 2$.  In most cases $\B(L)$ has valence 6, while $\C_i(L)$ and $C_x(L)$ each have valence 4.  The exception is when $q$ is divisible by 4 and $j = \frac{q}{2}-1$.   Here,  $\B(L)$ has valence 5 and $C_x(L)$ has valence 3.
%\end{theorem}
\begin{proof}
Because $j$ is odd, we know that the local corneration is equivalent to the standard odd local corneration.  Assuming, as we may, that $c_0\in L_v$,  we have that $L_v =\{c_{2i}|i\in \ZZ_{\frac{q}{2}}\}$.  The wedge-set for $\C_i$ is $\{(2i, 2i+1)|i\in \ZZ_{\frac{q}{2}}\}$, while that for $\C_x$ is $\{(2i-1, 2i)|i\in \ZZ_{\frac{q}{2}}\}$.

Thus, in $\C_i$, $(0, 1)$ connects $c_0$ to $c_{1-j}$ and $(j-1,j)$ connects $c_0$ to $c_{j-1}$.  In general, $c_i$ is connected to $c_{i+j-1}$, and to $c_{i-j+1}$.  Thus, each corner is connected by new edges to two other corners of $L$, giving it valence 4.  Further, a connected component in the local corneration consists of all $c_{i(j-1)}$, and this is all of $L_x$ if and only if $\gcd(q, j-1) = 2$, and so the graph is locally connected if  $\gcd(q, j-1) = 2$.  

 In $\C_x$, $(-1,0)$ connects $c_0$ to $c_{-1-j}$, and $ (j, j+1)$ connects $c_0$ to $c_{j+1}$, giving two new edges to $c_0$ unless $-1-j = j+1$. Thus, the valence of $\C_x(L)$ is 4 and that of $\B(L)$  is 6, unless $j = \frac{q}{2}-1$.  In that exceptional case, the valences are 3 and 5, respectively. Finally, since any $c_i$ is connected to $c_{i\pm (j+1)}$, $\C_x(L)$ is locally connected if and only if $\gcd(q, j+1) = 2$.
\end{proof}

Our final effort in this section is to determine valence and connectivity for graphs  $\B(L)$,$ C_i(L)$, $C_x(L)$ when $j$ is even.

\begin{theorem}\label{th:BCeven}
If $L$ is a transitive corneration of even width $2\le j<\frac{q}{2}$, then
\begin{enumerate}
\item $\B(L)$ is locally connected, 
\item $\C_i(L)$ is locally connected exactly when $\gcd(q, j-2) = 4$, 
\item $\C_x(L)$ is locally connected exactly when $\gcd(q, j+2) = 4$.  
\end{enumerate}
Furthermore:
\begin{enumerate}
\item[(4)]  $\B(L)$ has valence 6, unless $j=2$, in which case it has valence $5$,
\item[(5)] $\C_i(L)$ and $C_x(L)$ each have valence $4$ unless unless $j=2$, in which case it has valence $3$.
\end{enumerate}
\end{theorem}

%\begin{theorem}\label{th:BCeven}
%If $L$ is a transitive corneration of even width $2\le j<\frac{q}{2}$, then $\B(L)$ is locally connected, $\C_i(L)$ is locally connected exactly when $\gcd(q, j-2) = 4$, while $\C_x(L)$ is locally connected exactly when $\gcd(q, j+2) = 4$.  In most cases $\B(L)$ has valence 6, while $\C_i(L)$ and $C_x(L)$ each have valence 4.  The exception is when  $j = 2$.   Here,  $\B(L)$ has valence 5 and $C_i(L)$ has valence 3.
%\end{theorem}
\begin{proof}
Recall from the proof of Theorem \ref{pr:loccor}, that $j$ must be $2$ mod $4$, and the local corneration $L_v$ must be equivalent to the standard even local corneration.   We can assume, then, that $L_v$ is $\{c_i  | i {\rm (mod 4})\in \{0,3\}\}.$   Then the wedge-set for $\C_i(L)$ is all $(i,i+1)$ with $i$ (mod $4$) in $\{0,1,3\}$,  the wedge-set for $\C_x(L)$ is that for $\{1,2,3\}$, and, of course, the wedge-set for $\B(L)$ is the union of those, i.e., all wedges.  For $\C_i$, the wedge set includes $(-1, 0), (0, 1)$ and $(j-1, j)$, which connect $c_0$ to $c_{-1}, c_{1-j}, c_{-1}$ respectively.   Likewise, the wedge set includes $(-1, 0), (j-2, j-1)$ and $(j-1, j)$, which connect $c_{-1}$ to $c_{0}, c_{j-2}, c_{0}$ respectively.  Thus, each vertex has two new edges, giving it a valence of $4$, unless $j-2 = 0$, where the valence is 3..

	Further, each corner $(4i, 4i+j)$ is paired with $(4i-1, 4i+j-1)$, and paired vertices are joined by a new edge.  Moreover, one of each pair is joined to one from a pair $j-2$ places forward or backward around $v$; for instance, in the pair $\{c_{-1}, c_0\}$, $c_{-1}$ is joined to $c_{1-j}$, of the pair $\{c_{1-j}, c_{2-j}\}$, while $c_0$ is joined to $c_{j-2}$ of the pair $\{c_{j-3}, c_{j-2}\}$.   Thus, $\C_i(L)$ is locally connected if and only if $\gcd(q, j-2) = 4$.

Similar reasoning supports the conclusions abour $\C_x$ and $\B$.
\end{proof}

\subsection{Cubic splits}

One of our motivating examples is the construction of a cubic split graph from cycle structure.  To honor that heritage, we gather here the cases in which the constructed split graph has valence 3:
\begin{corollary}
Suppose that $L$ is a transitive $j$-corneration for $1 \leq j \leq \frac{q}{2}$.  Then:
\begin{itemize}
\item by Theorem \ref{th:Aval}, $\A(L)$ has valence $3$ if and only if $j = \frac{q}{2}$.
\item by Theorem \ref{th:straight}, $\B(L)$ has valence $3$ if and only if $j = \frac{q}{4}$.
\item by Theorem \ref{th:BCeven}, $\C_i(L)$ has valence $3$ if and only if $j = 2$.
\item by Theorem \ref{th:BCodd}, $\C_x(L)$ has valence $3$ if and only if  $q$ is divisible by $4$ and $j = \frac{q}{2}-1$.
\end{itemize}

\end{corollary}

\color{Black}

%%%%%%%%%%%%%%%%%%%%%%%%%%%%%%%%%%%%%%%%%%%%%%%%%%%%%%%%%%
%%%%%%%%%%   -------  Bibliography
%%%%%%%%%%%%%%%%%%%%%%%%%%%%%%%%%%%%%%%%%%%%%%%%%%%%%%%%%%

\end{document}